\documentclass[12pt,reqno,twoside]{amsart}%

\usepackage{amsmath,amsthm,amscd,amsthm,upref,indentfirst}
\usepackage{amsfonts,mathrsfs}
\usepackage{amssymb,amsbsy,bm}
\usepackage{graphicx}
\usepackage{amsaddr}
\usepackage{soul}

\usepackage[us,long,24hr]{datetime}

\usepackage{mathabx}
\usepackage[usenames,dvipsnames]{color}

\theoremstyle{plain}

\theoremstyle{definition}

\theoremstyle{remark}

\numberwithin{equation}{section}
\numberwithin{theorem}{section}

\theoremstyle{plain}
\newtheorem{Theorem}{Theorem}
\newtheorem{Assertion}[Theorem]{Assertion}

\newtheorem{Proposition}[Theorem]{Proposition}

\theoremstyle{definition}
\newtheorem{Definition}[Theorem]{Definition}

\newtheorem{Corollary}[Theorem]{Corollary}

\newtheorem{Notation}[Theorem]{Notation}

\theoremstyle{remark}
\newtheorem{Remark}[Theorem]{Remark}

\newtheorem{Example}[Theorem]{Example}

\numberwithin{equation}{section}
\numberwithin{Theorem}{section}

\normalfont\upshape

\usepackage{exscale}
\usepackage[scaled=0.86]{helvet}

\usepackage[scr=euler,scrscaled=1.1,bb=txof,bbscaled=1.1,cal=boondox,calscaled=1.1]{mathalfa}
\renewcommand{\mathbf}{\bm}
\renewcommand{\mathit}[1]{\mathscr #1}
\renewcommand{\mathfrak}[1]{{\textsc{\upshape #1}}}
\renewcommand{\mathtt}[1]{\scalebox{1.15}{\bf \texttt{\upshape#1}}}
\renewcommand{\emph}[1]{\mathtt{#1}}
\renewcommand{\mathrm}[1]{\scalebox{1.01}{\textbf{\upshape #1}}}

\usepackage{float}

\usepackage[justification=centering]{caption}

\numberwithin{equation}{section}
\numberwithin{theorem}{section}

\usepackage[normalem]{ulem}
\usepackage[square,comma]{natbib}

\usepackage{mparhack}

\usepackage{keyval}
\usepackage{totcount}
\newtotcounter{citnum} %From the package documentation
\def\oldbibitem{} \let\oldbibitem=\bibitem
\def\bibitem{\stepcounter{citnum}\oldbibitem}

\usepackage[verbose]{backref}
\backrefsetup{verbose=false}
\renewcommand*{\backref}[1]{}
\renewcommand*{\backrefalt}[4]{[{\tiny%
    \ifcase #1 \textsl{Not cited}%
          \or \textsl{Cited on page}~\textcolor{BrickRed}{#2}%
          \else \textsl{Cited on pages}~\textcolor{BrickRed}{#2}%
    \fi%
    }]}

\usepackage{epigraph}

\usepackage{fancyhdr}
\author{{\bfseries \small\scshape S\lowercase{teven} D\lowercase{uplij}}}

\address{% \small \scshape
Yantai Research Institute,
Harbin Engineering University,
265615 Yantai, China
\\and\\
Center for Information Technology,
University of M\"unster,
48149 M\"unster,
Germany
}
\email{\small \sf douplii@uni-muenster.de;
sduplij@gmail.com;
http://www.uni-muenster.de/IT.StepanDouplii}

\title{\large\bfseries\scshape
M\lowercase{ultiary gradings}}

\date{\textit{of start} December 9, 2025. \textit{Date}:
\textit{of completion} (ver.1) January 16, 2026; (ver. 2) March 9, 2026.
\newline
\mbox{}\hskip 1.16em
\textit{Total}:
36
references%,
}

\renewcommand{\refname}{\textsc{References}}

\let\origsection\section
\renewcommand{\section}[1]{\sectionmark{#1}\origsection{#1}}
\let\origsubsection\subsection
\renewcommand{\subsection}[1]{\subsectionmark{#1}\origsubsection{#1}}

\makeatletter
\renewenvironment{thebibliography}[1]{%
  \@xp\origsection\@xp*\@xp{\refname}%
  \normalfont\footnotesize\labelsep .9em\relax
  \renewcommand\theenumiv{\arabic{enumiv}}\let\p@enumiv\@empty
  \vspace*{-5pt}% NEW
  \list{\@biblabel{\theenumiv}}{\settowidth\labelwidth{\@biblabel{#1}}%
    \leftmargin\labelwidth \advance\leftmargin\labelsep
    \usecounter{enumiv}}%
  \sloppy \clubpenalty\@M \widowpenalty\clubpenalty
  \sfcode`\.=\@m
}{%
  \def\@noitemerr{\@latex@warning{Empty `thebibliography' environment}}%
  \endlist
}
\makeatother

\subjclass[2010]{13A02, 13A30, 16W50, 17A70,
 20H20, 20N10, 20N15, 81Q60}
\keywords{grading, graded algebra, arity, polyadic ring, $n$-ary group, polyadic algebra, querelement}

\begin{document}
\mbox{}
\vskip 1.8cm
\begin{abstract}
%TCIDATA{OutputFilter=latex2.dll}
%TCIDATA{Version=5.50.0.2953}
%TCIDATA{LaTeXparent=0,0,example.TEX}

This article develops a comprehensive theory of multiary graded polyadic algebras, extending the classical concept of group-graded algebras to higher-arity structures. We introduce the notion of grading by multiary groups and investigate various compatibility conditions between the arity of algebra operations and grading group operations. Key results include quantization rules connecting arities, classification of graded homomorphisms, the First Isomorphism Theorem for graded polyadic algebras and concrete examples including ternary superalgebras and polynomial algebras over $n$-ary matrices. The theory reveals fundamentally new phenomena not present in the binary case, such as the existence of higher power gradings and nontrivial constraints on arity compatibility.

\end{abstract}

\maketitle

\thispagestyle{empty}
\mbox{}
\vspace{-0.5cm}
%\begin{small}
\tableofcontents
%\end{small}%
\newpage

\pagestyle{fancy}

\addtolength{\footskip}{15pt}

\renewcommand{\sectionmark}[1]{%
\markboth{
{ \scshape #1}}{}}

\renewcommand{\subsectionmark}[1]{%
\markright{
\mbox{\;}\\[5pt]
\textmd{#1}}{}}

\fancyhead{}
\fancyhead[EL,OR]{\leftmark}
\fancyhead[ER,OL]{\rightmark}
\fancyfoot[C]{\scshape -- \textcolor{BrickRed}{\thepage} --}
\fancyfoot[R]{}
%\fancyfoot[L]{{\small \hours}}

\renewcommand\headrulewidth{0.5pt}
\fancypagestyle {plain1}{ %
\fancyhf{}
\renewcommand {\headrulewidth }{0pt}
\renewcommand {\footrulewidth }{0pt}
}

\fancypagestyle{plain}{ %
\fancyhf{}
\fancyhead[C]{\scshape S\lowercase{teven} D\lowercase{uplij} \hskip 0.7cm \MakeUppercase{Polyadic Hopf algebras and quantum groups}}
\fancyfoot[C]{\scshape - \thepage  -}
\renewcommand {\headrulewidth }{0pt}
\renewcommand {\footrulewidth }{0pt}
}

\fancypagestyle{fancyref}{ %
\fancyhf{} % remove everything
\fancyhead[C]{\scshape R\lowercase{eferences} }
\fancyfoot[C]{\scshape -- \textcolor{BrickRed}{\thepage} --}
\renewcommand {\headrulewidth }{0.5pt}
\renewcommand {\footrulewidth }{0pt}
}

\fancypagestyle{emptyf}{
\fancyhead{}
\fancyfoot[C]{\scshape -- \textcolor{BrickRed}{\thepage} --}
\renewcommand{\headrulewidth}{0pt}
}
\thispagestyle{emptyf}
\section{\textsc{Introduction}}

The theory of graded algebras and rings~\cite{nas/oys,nas/oys04,hazrat}
represents a fundamental chapter in modern algebra~\cite{lang,rotman}, with
profound applications across mathematics~\cite{bah/seh/zai,kelarev,haz/gub,lam1991} and theoretical physics~\cite{wes/bag,del/eti/fre,terning,gre/sch/wit,kaku3}. Classically, a~grading
on an algebra $\mathcal{A}$ over a field $\Bbbk$ by a group $\mathsf{G}$
consists of a decomposition $\mathcal{A}=\bigoplus_{\mathsf{g}\in\mathsf{G}%
}\mathcal{A}(\mathsf{g})$ into homogeneous components such that the
multiplication respects the group structure $\mathcal{A}(\mathsf{g}%
)\cdot\mathcal{A}(\mathsf{h})\subseteq\mathcal{A}(\mathsf{g}+\mathsf{h})$.
This elegant framework encompasses numerous important structures including
superalgebras (graded   by $\mathbb{Z}_{2}$) \cite{berezin,kac3,leites13},
polynomial algebras (graded by degree) \cite{mat59,russell}, group algebras
and rings~\cite{zal/mik,passman,bovdi}.

Recent developments in polyadic (multiplace) algebra~\cite{duplij2022} have
opened new avenues for generalization, where binary operations are replaced by
$n$-ary ones. Building upon our previous work on the polyadic generalization of
group rings $\mathcal{R}\left[  \mathsf{G}\right]  $ \cite{dup2025h}, we now
extend the related concept of the grading to the polyadic realm. Polyadic
structures exhibit rich and often surprising behavior: polyadic groups may
lack identity elements, polyadic fields may have no zero or unit, and
associativity takes more intricate forms. These features necessitate a
fundamental revision of many standard algebraic~constructions.

In this work, we introduce and systematically develop the theory of multiary
graded polyadic algebras. Our approach follows the ``arity freedom principle''
\cite{duplij2022}, allowing initial arities to be arbitrary, with~structural
constraints emerging from compatibility conditions. The~main contributions are as follows:

\begin{enumerate}
\item [(1)] A general definition of multiary $\mathsf{G}$-graded polyadic $\Bbbk
$-algebras, where both the algebra operations and grading group operations may
have arbitrary arities ({Section} \ref{sec-gen}).

\item [(2)] Quantization rules connecting the arities of algebra multiplication
($n$), grading group multiplication ($n^{\prime}$), and~algebra addition ($m$)
with the order of the grading group ({{Theorems} 
} \ref{theor-G=l} and
\ref{theor-l=l}).

\item [(3)] A classification of graded homomorphisms between multiary graded
algebras, and~the First Isomorphism Theorem for graded polyadic algebras is
proved ({Section} \ref{sec-hom}).

\item [(4)] Concrete examples including

\begin{itemize}
\item Derived ternary superalgebras with binary $\mathbb{Z}_{2}$-grading
({Section} \ref{sec-super});

\item Strictly nonderived ternary superalgebras with ternary grading groups ({Example} \ref{ex-tern});

\item Polynomial algebras over $n$-ary matrices graded by polyadic integers
$\mathbb{Z}^{[m^{\prime\prime},n^{\prime\prime}]}$ ({Section}
\ref{sec-poly});

\item Higher power multiary gradings with $n^{\prime}\neq n$ ({Section}
\ref{sec-high}).
\end{itemize}

\item [(5)] Investigation of support properties, strong grading conditions, and~the
relationship between graded components and polyadic powers.
\end{enumerate}

Our results reveal that the transition from binary to polyadic grading
introduces qualitatively new phenomena. For~instance, the~compatibility
condition between algebra $n$-ary multiplication and grading group $n^{\prime
}$-ary multiplication leads to ``quantization'' rules like $n^{\prime}=n$ for
strongly graded algebras, while higher power gradings satisfy $\ell
_{n^{\prime}}(n^{\prime}-1)=\ell_{n}(n-1)$, where $\ell$ denotes polyadic
power. These constraints have no analogs in classical graded algebra~theory.

The article is structured as follows: {Section} \ref{sec-prelim}
reviews the essential polyadic notations and preliminaries. {Section}
\ref{sec-gen} presents the general theory of multiary graded polyadic
algebras. {Section} \ref{sec-hom} develops the theory of graded
homomorphisms. {Sections} \ref{sec-super}--\ref{sec-high} provide detailed examples and applications, demonstrating the
richness and versatility of the~theory.

\section{\textsc{Preliminaries}\label{sec-prelim}}

A grading on an algebra $\mathcal{A}$ (over a field $\Bbbk$, say) by a group
$\mathsf{G}$ (called the grading group) is a way to decompose $\mathcal{A}$
into \textquotedblleft layers\textquotedblright\ indexed by elements of
$\mathsf{G}$. Formally, the~algebra $\mathcal{A}$ is $\mathsf{G}$-graded, if~

\begin{enumerate}
\item [(1)] $\mathcal{A}$ {is a} 
 direct sum of $\Bbbk$-vector spaces%
\begin{equation}
\mathcal{A}=\bigoplus_{\mathsf{g}\in\mathsf{G}}\mathcal{A}\left(
\mathsf{g}\right)  , \label{ab}%
\end{equation}
where each $\mathcal{A}\left(  \mathsf{g}\right)  $ is called the homogeneous
component of degree $\mathsf{g}$, and~they do not intersect: $\mathcal{A}%
\left(  \mathsf{g}\right)  \cap\mathcal{A}\left(  \mathsf{h}\right)
=\varnothing$, if~$\mathsf{g}\neq\mathsf{h}$. So, every element of
$\mathcal{A}$ can be uniquely written {as} 
 $\mathbf{a}=\sum_{\mathsf{g}%
\in\mathsf{G}}\mathbf{a}\left(  \mathsf{g}\right)  $, with~$\mathbf{a}\left(
\mathsf{g}\right)  \in\mathcal{A}\left(  \mathsf{g}\right)  $, notation
$\mathsf{g}=\deg\left(  \mathbf{a}\right)  $, and~only finitely many nonzero
elements $\mathbf{a}\left(  \mathsf{g}\right)  \neq\mathbf{0}$, such that, in
each component, it has finite~support.

\item [(2)] The multiplication in $\mathcal{A}$ respects the grading, that is, for
product of algebra components
\begin{equation}
\mathcal{A}\left(  \mathsf{g}\right)  \cdot\mathcal{A}\left(  \mathsf{h}%
\right)  \subseteq\mathcal{A}\left(  \mathsf{g}+\mathsf{h}\right)
,\ \ \ \mathsf{g},\mathsf{h}\in\mathsf{G}, \label{ag}%
\end{equation}
where the grading group $\mathsf{G}$ is typically abelian ($\mathbb{Z}%
,\mathbb{Z}_{n},\mathbb{Z}^{\times n}$), written additively, with~the identity
$\mathsf{0}$. Thus, $\mathcal{A}\left(  \mathsf{0}\right)  \cdot
\mathcal{A}\left(  \mathsf{g}\right)  \subseteq\mathcal{A}\left(
\mathsf{g}\right)$, and $\mathcal{A}\left(  \mathsf{g}\right)  \cdot
\mathcal{A}\left(  \mathsf{0}\right)  \subseteq\mathcal{A}\left(
\mathsf{g}\right)  $. The~grading is strong if,~instead of (\ref{ag}), the
equality $\mathcal{A}\left(  \mathsf{g}\right)  \cdot\mathcal{A}\left(
\mathsf{h}\right)  =\mathcal{A}\left(  \mathsf{g}+\mathsf{h}\right)  $ holds
(used for crossed product algebras).
\end{enumerate}

The grading is \textquotedblleft natural\textquotedblright, if~algebra
homomorphisms preserve grading, that is, for any algebras map $\Phi
:\mathcal{A}\rightarrow\mathcal{B}$, the~homogeneity is preserved $\Phi\left(
\mathcal{A}\left(  \mathsf{g}\right)  \right)  \subseteq\mathcal{B}\left(
\mathsf{g}\right)  $, for~any $\mathsf{g}\in\mathsf{G}$. If~the algebra is
associative, we can extend the associativity across grades to preserve degrees. If~
the algebra element is invertible, then $\deg\left(  \mathbf{a}^{-1}\right)
=-\mathsf{g}$ (in additive notation).

The most common examples of the graded algebra are polynomial algebra
$\Bbbk\left[  x\right]  $, exterior algebra $\wedge V$ and superalgebra and group
algebra $\Bbbk\left[  \mathsf{G}\right]  $ (any $\Bbbk\left[  \mathsf{G}%
\right]  $ is a $\mathsf{G}$-graded algebra, but~not vice-versa).

To generalize the concept of graded algebra to the polyadic case, we need some
polyadic notation~\cite{duplij2022}, to~be self-consistent~here.

\begin{Notation}
An algebraic structure $\mathcal{M}$ will be denoted by%
\begin{equation}
\mathcal{M}=\mathcal{M}_{order}^{\left[  arity\right]  }\left(
parameter\right)  ,
\end{equation}
where $arity$ represents the number of operation places under consideration,
that is, arity(ies), $order$ denotes the number of the carrier element or power of
underlying set(s), $parameter$ gives additional variable(s), which characterize
the structure $\mathcal{M}$, and~they are used when necessary.
\end{Notation}

Let $A$ be a set and $\mu^{\left[  n\right]  }:A^{\times n}\rightarrow A$ be a
polyadic multiplication. An~$n$-ary magma is $A$ closed under $\mu^{\left[
n\right]  }$ and denoted by $\mathcal{M}^{\left[  n\right]  }=\left\langle
A\mid\mu^{\left[  n\right]  }\right\rangle $. The~admissible word length in
$\mathcal{M}^{\left[  n\right]  }$ is \textquotedblleft
quantized\textquotedblright, which means%
\begin{equation}
w=\ell\left(  n-1\right)  +1, \label{w}%
\end{equation}
where $\ell$ is a polyadic power, that is, the \textquotedblleft number of
multiplications\textquotedblright. For~any element $a\in A$, the polyadic power
$\ell$ is denoted by%
\begin{equation}
x^{\left\langle \ell\right\rangle }=\mu^{\left[  n\right]  \circ\ell}\left[
\overset{\ell\left(  n-1\right)  +1}{\overbrace{a,a,\ldots,a}}\right]  .
\label{xl}%
\end{equation}

An $n$-ary operation is strictly nonderived, if~it cannot be composed from a
binary operation without fixing elements or additional conditions. A~polyadic
identity $e\in A$ is defined as (if it exists)%
\begin{equation}
\mu^{\left[  n\right]  }\left[  \overset{n-1}{\overbrace{e,e,\ldots,e}%
},a\right]  =a, \label{ea}%
\end{equation}
and the polyadic zero $z\in A$ is%
\begin{equation}
\mu^{\left[  n\right]  }\left[  \overset{n-1}{\overbrace{a,a,\ldots,a}%
},z\right]  =z. \label{az}%
\end{equation}

A polyadic semigroup in a one-set one-operation {structure} %MDPI: Please check if all variables are unified with the format of the equations in the paper.
%%%%%%Author's ANSWER: Checked - OK.
 $\mathcal{S}%
^{\left[  n\right]  }=\left\langle A\mid\mu^{\left[  n\right]  }%
,assoc\right\rangle$, which is totally associative%
\begin{equation}
\mu^{\left[  n\right]  \circ2}\left[  \mathbf{a},\mathbf{b},\mathbf{c}\right]
=\mu^{\left[  n\right]  }\left[  \mathbf{a},\mu^{\left[  n\right]  }\left[
\mathbf{b}\right]  ,\mathbf{c}\right]  =invariant, \label{mab}%
\end{equation}
with respect to the position of the middle multiplication, where
$\mathbf{a},\mathbf{b},\mathbf{c}$ are polyads (sequences of elements) of the
needed sizes. Note, that in higher arity case the associativity can be partial~\cite{sok1}, which can lead to more complicated structural consequences  (see,
e.g., \cite{kum/wat,sam/geo}).

The querelement $\bar{a}$ for $a\in S$ is defined by%
\begin{equation}
\mu^{\left[  n\right]  }\left[  \overset{n-1}{\overbrace{a,a,\ldots,a}}%
,\bar{a}\right]  =a, \label{aq}%
\end{equation}
which can be treated as an unary queroperation $\bar{\mu}^{\left[  1\right]
}\left[  a\right]  =\overline{\left(  a\right)  }=\bar{a}$. If~every element
of an $n$-ary semigroup $\mathcal{S}$ has the unique querelement, it is called a
polyadic ($n$-ary) group $\mathcal{G}^{\left[  n\right]  }=\left\langle
A\mid\mu^{\left[  n\right]  },\overline{\left(  {}\right)  }%
,assoc\right\rangle $. A~polyadic or $\left(  m,n\right)  $-ring is a one-set
two operation algebraic structure $\mathcal{R}^{\left[  m,n\right]
}=\left\langle A\mid\nu^{\left[  m\right]  },\mu^{\left[  n\right]
},assoc,distr\right\rangle $, such that $m$-ary addition $\nu^{\left[
m\right]  }$ and $n$-ary multiplication $\mu^{\left[  n\right]  }$ satisfy
distributivity~\cite{lee/but}.

\begin{Example}
[Ternary polynomial graded algebra]\label{ex-tern-pol}A natural example of a
ternary $\mathbb{Z}$-graded algebra is a polynomial ring of three variables
$\mathcal{A}^{\left[  2,3\right]  }=\Bbbk\left[  x,y,z\right]  $ over a field
$\Bbbk=\mathbb{R}$, $x,y,z\in\mathbb{R}$, and the~addition and multiplication by
scalar are binary for the algebra $\mathcal{A}^{\left[  2,3\right]  }$. The~
grading of elements is defined by $\deg\left(  1\right)  =0\in\mathbb{Z}$,
$\deg\left(  x\right)  =1\in\mathbb{Z}$, and~the monomial degree is%
\begin{equation}
\deg\left(  x^{p^{\prime}}y^{q^{\prime}}z^{r^{\prime}}\right)  =p^{\prime
}+q^{\prime}+r^{\prime}, \label{pqr}%
\end{equation}
where $p^{\prime},q^{\prime},r^{\prime}\in\mathbb{Z}$, and~therefore, the
grading group is binary $\mathsf{G}=\mathbb{Z}$. The~direct sum decomposition
(\ref{ab}) becomes%
\begin{equation}
\mathcal{A}^{\left[  2,3\right]  }=\bigoplus_{p\in\mathbb{Z}}\mathcal{A}%
\left(  p\right)  , \label{a3p}%
\end{equation}
where $\mathcal{A}^{\left[  3\right]  }\left(  p\right)  $ is the homogeneous
component of degree $p\in\mathbb{Z}$, and~obviously, $\mathcal{A}\left(
p\right)  \cap\mathcal{A}\left(  q\right)  =\varnothing$ for $p\neq
q\in\mathbb{Z}$. The~ternary product in algebra $\mathcal{A}^{\left[
3\right]  }$ is defined as the ordinary multiplication $\left(  \cdot\right)
$ of polynomials in three variables $f\left(  x,y,z\right)  \in\Bbbk\left[
x,y,z\right]  $ as%
\begin{equation}
\mu^{\left[  3\right]  }\left[  f,g,h\right]  =f\cdot g\cdot h\in
\mathcal{A}^{\left[  2,3\right]  },\ \ \ f,g,h\in\mathcal{A}^{\left[
2,3\right]  }, \label{m3f}%
\end{equation}
such that $\mu^{\left[  3\right]  }$ is trilinear and respects grading,
because from (\ref{pqr}) and (\ref{m3f}), cf. the binary case (\ref{ag}),   it
follows that%
\begin{equation}
\mu^{\left[  3\right]  }\left[  \mathcal{A}\left(  p\right)  ,\mathcal{A}%
\left(  q\right)  ,\mathcal{A}\left(  r\right)  \right]  =\mathcal{A}\left(
p+q+r\right)  . \label{m3a}%
\end{equation}

Because there is equality in (\ref{m3a}), the~grading is strong. Thus, the
ternary graded algebra $\mathcal{A}^{\left[  2,3\right]  }$ is commutative and
ternary associative, it is non-truncated because there are no relations, and the~
polynomial ring is free. It is infinite-dimensional as a vector space, but
each homogeneous component $\mathcal{A}\left(  p\right)  $ of the direct sum
decomposition (\ref{a3p}) is finite-dimensional,   and~$\dim\mathcal{A}\left(
p\right)  =\binom{p+2}{2}$.
\end{Example}

\section{\textsc{General Grading of Polyadic Algebra by
Multiary Group}\label{sec-gen}}

Let us generalize the grading concept to the polyadic case by consideration of
higher arity in both the algebra and grading group. According to the
\textquotedblleft arity freedom principle\textquotedblright\ \cite{duplij2022},
the initial arities can be arbitrary; then, the structural constraints appear
from the general dependencies leading to \textquotedblleft quantization
rules\textquotedblright. Polyadic structures can display unusual properties
and behaviors. For~example, some polyadic groups have no identity element at
all, while others may have multiple identities. Likewise, there are polyadic
fields that lack a zero, a~unit, or~both~\cite{duplij2022}. These features can
require a substantial rethinking of many standard mathematical~statements.

A polyadic algebra over a polyadic field is the 2-set 5-higher arity operation
algebraic structure satisfying compatibility axioms~\cite{dup2019,duplij2022}.
However, in~our consideration here, we will use the binary field $\Bbbk$ and
consider only the corresponding $\Bbbk$-algebra $m$-ary addition $\mathbf{\nu
}_{a}^{\left[  m\right]  }:A^{\times m}\rightarrow A$ and $n$-ary
multiplication $\mathbf{\mu}_{a}^{\left[  n\right]  }:A^{\times n}\rightarrow
A$, and~therefore, for $n$-ary algebra we denote%
\begin{equation}
\mathcal{A}^{\left[  m,n\right]  }=\left\langle A\mid\mathbf{\nu}_{a}^{\left[
m\right]  },\mathbf{\mu}_{a}^{\left[  n\right]  }\right\rangle , \label{an}%
\end{equation}
where $A$ is its underlying set, and~zero or unit are not necessary. The~
grading group now becomes the $n^{\prime}$-ary group%
\begin{equation}
\mathsf{G}^{\left[  n^{\prime}\right]  }=\left\langle G\mid\mu_{g}^{\left[
n^{\prime}\right]  },\overline{\left(  \ \right)  }\right\rangle , \label{g}%
\end{equation}
where $G$ is the group underlying set, the~$n^{\prime}$-ary multiplication is
$\mu_{g}^{\left[  n^{\prime}\right]  }:G^{\times n^{\prime}}\rightarrow G$, and
$\overline{\left(  \ \right)  }$ is the unary operation of taking querelement
(polyadic inverse (\ref{aq})), while no identity element (\ref{ea}) is
necessary but~can exist for some polyadic groups; see, e.g.~\cite{duplij2022}%
. Similarly to the binary case, we assume that the grading $n^{\prime}$-ary
group $\mathsf{G}^{\left[  n^{\prime}\right]  }$ is~commutative.

\begin{Notation}
We call the generic algebraic structure with multiplace operations (after~\cite{pos}) \textquotedblleft polyadic\textquotedblright\ (ring, field,
algebra, group), while for the grading group, just for distinctiveness,
convenience and clarity, we use its Latin synonym \textquotedblleft
multiary\textquotedblright\ \cite{kas}.
\end{Notation}

\begin{Remark}
\label{rem-nonder}To be far from the well-known binary case, we assume both
{multiplications}  $\mathbf{\mu}_{a}^{\left[  n\right]  }$ and $\mu_{g}^{\left[
n^{\prime}\right]  }$ are nonderived (or not strictly derived), i.e.,~they are
not consequent compositions of binary operations only (the weaker $b$-derived
multiplication is a different story~\cite{dor3}).
\end{Remark}

The expansion into the direct sum for the polyadic algebra $\mathcal{A}%
^{\left[  m,n\right]  }$ (\ref{an}) does not differ from that of the binary
algebras (\ref{ab}). The~crucial peculiarities come from the polyadic analog
of the compatibility condition (\ref{ag}), which gives several kinds of the
polyadic graded~algebras.

\begin{Definition}
A multiary $\mathsf{G}$-graded polyadic $\Bbbk$-algebra is the direct sum
decomposition of   $\Bbbk$-vector spaces%
\begin{equation}
\mathcal{A}^{\left[  m,n\right]  }=\bigoplus_{\mathsf{g}_{i}\in\mathsf{G}%
^{\left[  n^{\prime}\right]  }}\mathcal{A}\left(  \mathsf{g}_{i}\right)  ,
\label{a}%
\end{equation}
such that the $n$-ary multiplication in the algebra respects the $n^{\prime}%
$-ary multiplication in the   grading group%
\begin{equation}
\mathbf{\mu}_{a}^{\left[  n\right]  }\left[  \mathcal{A}\left(  \mathsf{g}%
_{1}\right)  ,\mathcal{A}\left(  \mathsf{g}_{2}\right)  ,\ldots,\mathcal{A}%
\left(  \mathsf{g}_{n}\right)  \right]  \subseteq\mathcal{A}\left(  \mu
_{g}^{\left[  n^{\prime}\right]  }\left[  \mathsf{g}_{1},\mathsf{g}_{2}%
,\ldots,\mathsf{g}_{n^{\prime}}\right]  \right)  , \label{ma}%
\end{equation}
where $\mathcal{A}\left(  \mathsf{g}_{i}\right)  $ is the $i$th component of
the decomposition (\ref{a}).

A polyadic algebra is strongly graded, if, as in (\ref{ma}), the equality is%Please confirm meaning is retained.
\begin{equation}
\mathbf{\mu}_{a}^{\left[  n\right]  }\left[  \mathcal{A}\left(  \mathsf{g}%
_{1}\right)  ,\mathcal{A}\left(  \mathsf{g}_{2}\right)  ,\ldots,\mathcal{A}%
\left(  \mathsf{g}_{n}\right)  \right]  =\mathcal{A}\left(  \mu_{g}^{\left[
n^{\prime}\right]  }\left[  \mathsf{g}_{1},\mathsf{g}_{2},\ldots
,\mathsf{g}_{n^{\prime}}\right]  \right)  . \label{mas}%
\end{equation}

\end{Definition}

\begin{Definition}
The support of a polyadic $\mathsf{G}$-graded algebra \textsl{supp}$\left(
\mathcal{A}^{\left[  m,n\right]  }\right)  $ is the set of all grading group
elements $\mathsf{g}_{i}\in\mathsf{G}^{\left[  n^{\prime}\right]  }$ such that
the corresponding component $\mathcal{A}\left(  \mathsf{g}_{i}\right)  $ as
$\Bbbk$-vector space contributes to the direct sum decomposition
(\ref{a}).
\end{Definition}

The number of \textquotedblleft non-zero\textquotedblright\ (contributing)
summands in the decomposition (\ref{a}) is exactly the cardinality of the
support; that is, $\left\vert \text{\textsl{supp}}\left(  \mathcal{A}^{\left[
m,n\right]  }\right)  \right\vert $. Because~\textsl{supp}$\left(
\mathcal{A}^{\left[  m,n\right]  }\right)  \subseteq\mathsf{G}^{\left[
n^{\prime}\right]  }$, the~number of summands in (\ref{a}) is $\left\vert
\text{\textsl{supp}}\left(  \mathcal{A}^{\left[  m,n\right]  }\right)
\right\vert \leq\left\vert G\right\vert $.

\begin{Assertion}
For strongly $\mathsf{G}$-graded polyadic algebra (\ref{mas}),%
\begin{equation}
\left\vert \text{\textsl{supp}}\left(  \mathcal{A}^{\left[  m,n\right]
}\right)  \right\vert =\left\vert G\right\vert . \label{sg}%
\end{equation}

\end{Assertion}

\begin{proof}
If $\mathcal{A}^{\left[  m,n\right]  }$ is strongly graded and
\textquotedblleft non-zero\textquotedblright, then every $\mathcal{A}\left(
\mathsf{g}_{i}\right)  $ is \textquotedblleft non-zero\textquotedblright, and~
therefore, the number of summands coincides with the graded group order
$\left\vert G\right\vert $.
\end{proof}

\begin{Proposition}
\label{prop-n=n}The strong polyadic algebra arity of multiplication and the
arity of multiary grading group coincide:%
\begin{equation}
n^{\prime}=n. \label{mn}%
\end{equation}

\end{Proposition}

\begin{proof}
This follows from the construction (\ref{mas}) by equating the number of components
as the multiplier in the algebra product $\mathbf{\mu}_{a}^{\left[  n\right]  }$
and number of the elements in the grading group   product $\mu_{g}^{\left[
n^{\prime}\right]  }$.
\end{proof}

\begin{Remark}
\label{rem-assoc}If the grading group has the identity $\mathsf{e}%
\in\mathsf{G}^{\left[  n^{\prime}\right]  }$, at~first glance, we can obtain
the unequally possibility $n^{\prime}>n$ by adding $\mathsf{e}$ to the
remaining $n^{\prime}-n$ places in (\ref{mas}) to get $\mu_{g}^{\left[
n^{\prime}\right]  }\left[  \mathsf{g}_{1},\mathsf{g}_{2},\ldots
,\mathsf{g}_{n},\overset{n^{\prime}-n}{\overbrace{\mathsf{e},\ldots
,\mathsf{e}}}\right]  $, but~this destroys the coincidence of the total polyadic
associativity in both sides of (\ref{mas}).
\end{Remark}

\begin{Remark}
\label{rem-hom}The strong consistency condition (\ref{mas}) can be interpreted
such that the components $\mathcal{A}\left(  \mathsf{g}_{i}\right)  $ for
$n^{\prime}=n$ satisfy the polyadic homomorphism relation~\cite{duplij2022}.
\end{Remark}

In the next section, it will be shown that there exist such multiary
gradings, for~which the arities $n^{\prime}$ and $n$ are different, which can
lead to unusual structural consequences, which are  impossible in the ordinary
binary~case.

In the manifest form, the~direct sum decomposition (\ref{a}) means that each
element of $\mathcal{A}^{\left[  m,n\right]  }$ can be written as a finite sum
(using the algebra $m$-ary addition) of   homogeneous elements%
\begin{equation}
\mathbf{a}=\mathbf{\nu}_{a}^{\left[  m\right]  }\left[  \mathbf{a}\left(
\mathsf{g}_{1}\right)  ,\mathbf{a}\left(  \mathsf{g}_{2}\right)
,\ldots,\mathbf{a}\left(  \mathsf{g}_{m}\right)  \right]  ,\ \ \mathbf{a}%
\in\mathcal{A}^{\left[  m,n\right]  },\ \mathbf{a}\left(  \mathsf{g}%
_{i}\right)  \in\mathcal{A}\left(  \mathsf{g}_{i}\right)  ,\ \mathsf{g}_{i}%
\in\mathsf{G}^{\left[  n^{\prime}\right]  }. \label{am}%
\end{equation}

However, this expression contains only one algebra addition operation
$\mathbf{\nu}_{a}^{\left[  m\right]  }$, which corresponds to the one
summation in the binary case as $\mathbf{a}=\mathbf{a}\left(  \mathsf{g}%
_{1}\right)  +\mathbf{a}\left(  \mathsf{g}_{2}\right)  $, while the binary
decomposition (\ref{ab}) contains $\left\vert \text{\textsl{supp}}\left(
\mathcal{A}^{\left[  m,n\right]  }\right)  \right\vert -1$ summations. Denote
by $\ell_{m}$ the polyadic power (\ref{xl}) of $m$-ary addition in the
polyadic algebra $\mathcal{A}^{\left[  m,n\right]  }$. Then, instead of
(\ref{am}), we have the decomposition of each element as%
\begin{equation}
\mathbf{a}=\mathbf{\nu}_{a}^{\left[  m\right]  \circ\ell_{m}}\left[
\mathbf{a}\left(  \mathsf{g}_{1}\right)  ,\mathbf{a}\left(  \mathsf{g}%
_{2}\right)  ,\ldots,\mathbf{a}\left(  \mathsf{g}_{\ell_{m}\left(  m-1\right)
+1}\right)  \right]  ,\ \ \mathbf{a}\in\mathcal{A}^{\left[  m,n\right]
},\ \mathbf{a}\left(  \mathsf{g}_{i}\right)  \in\mathcal{A}\left(
\mathsf{g}_{i}\right)  ,\ \mathsf{g}_{i}\in\mathsf{G}^{\left[  n^{\prime
}\right]  }. \label{al}%
\end{equation}

For a strongly graded algebra$\mathcal{\ A}^{\left[  m,n\right]  }$, there are
no zeroes among homogeneous components $\mathbf{a}\left(  \mathsf{g}%
_{i}\right)  $. In~this case, the~support is not just a subgroup; it is the
entire group \textsl{supp}$\left(  \mathcal{A}^{\left[  m,n\right]  }\right)
=\mathsf{G}^{\left[  n^{\prime}\right]  }$.

\begin{Theorem}
\label{theor-G=l}In the polyadic strongly $\mathsf{G}$-graded algebras, the
order of the finite grading group $\left\vert G\right\vert $ and the arity of
algebra addition $m$ are connected by%
\begin{equation}
\left\vert G\right\vert =\ell_{m}\left(  m-1\right)  +1. \label{sl}%
\end{equation}

\end{Theorem}

\begin{proof}
This follows from the admissible length $\ell_{m}\left(  m-1\right)  +1$ of a
word for $m$-ary operation repeated $\ell_{m}$ times (\ref{w}) and from
(\ref{sg}), (\ref{al}).
\end{proof}

In the case of the binary algebra addition $m=2$, the~number of summands in
the decomposition becomes $\ell_{m}+1$, where $\ell_{m}$ is the number of sum
signs (\ref{xl}). In~most cases, \textsl{supp}$\left(  \mathcal{A}^{\left[
m,n\right]  }\right)  $ is a subgroup of $\mathsf{G}^{\left[  n^{\prime
}\right]  }\mathsf{.}$ If $\mathcal{A}^{\left[  m,n\right]  }\mathcal{\ }$is
unital with the unity $\mathbf{e}$, then $\mathbf{e}\in\mathcal{A}\left(
\mathsf{e}\right)  $, where $\mathsf{e}$ is the identity of the group
$\mathsf{G}^{\left[  n^{\prime}\right]  }$ (if it exists), and then $\mathsf{e}\in
$\textsl{supp}$\left(  \mathcal{A}^{\left[  m,n\right]  }\right)  $.

\section{\textsc{Polyadic Graded~Homomorphisms}\label{sec-hom}}

Consider another multiary $\mathsf{H}$-graded polyadic algebra,%
\begin{equation}
\mathcal{B}^{\left[  m,n\right]  }=\bigoplus_{\mathsf{h}_{i}\in\mathsf{H}%
^{\left[  n^{\prime}\right]  }}\mathcal{B}\left(  \mathsf{h}_{i}\right)  ,
\label{b}%
\end{equation}
which has the same arity shape as $\mathcal{A}^{\left[  m,n\right]  }$
(\ref{an}).

\begin{Definition}
A polyadic graded homomorphism is a pair $\dbinom{\Phi}{\Psi}$ of maps
$\Phi:\mathcal{A}^{\left[  m,n\right]  }\rightarrow\mathcal{B}^{\left[
m,n\right]  }$ and $\Psi:\mathsf{G}^{\left[  n^{\prime}\right]  }%
\mathsf{\rightarrow H}^{\left[  n^{\prime}\right]  }$, which respect $m$-ary
additions%
\begin{equation}
\Phi\left(  \mathbf{\nu}_{a}^{\left[  m\right]  }\left[  \mathbf{a}%
_{1},\mathbf{a}_{2},\ldots,\mathbf{a}_{m}\right]  \right)  =\mathbf{\nu}%
_{b}^{\left[  m\right]  }\left[  \Phi\left(  \mathbf{a}_{1}\right)
,\Phi\left(  \mathbf{a}_{2}\right)  ,\ldots,\Phi\left(  \mathbf{a}_{m}\right)
\right]  , \label{fn}%
\end{equation}
and $n$-ary multiplications%
%\vspace{-10pt}
%\begin{adjustwidth}{-\extralength}{0cm}
\begin{equation}
\Phi\left(  \mathbf{\mu}_{a}^{\left[  n\right]  }\left[  \mathbf{a}%
_{1},\mathbf{a}_{2},\ldots,\mathbf{a}_{n}\right]  \right)  =\mathbf{\mu}%
_{b}^{\left[  n\right]  }\left[  \Phi\left(  \mathbf{a}_{1}\right)
,\Phi\left(  \mathbf{a}_{2}\right)  ,\ldots,\Phi\left(  \mathbf{a}_{n}\right)
\right]  ,\ \ \ \mathbf{a}_{j}\in\mathcal{A}^{\left[  m,n\right]  }%
,\ \Phi\left(  \mathbf{a}_{j}\right)  \in\mathcal{B}^{\left[  m,n\right]  },
\label{fm}%
\end{equation}
%\end{adjustwidth}
in algebras and $n^{\prime}$-ary multiplications in the grading groups%
\begin{equation}
\Psi\left(  \mathbf{\mu}_{g}^{\left[  n^{\prime}\right]  }\left[
\mathsf{g}_{1},\mathsf{g}_{2},\ldots,\mathsf{g}_{n^{\prime}}\right]  \right)
=\mu_{g}^{\left[  n^{\prime}\right]  }\left[  \Psi\left(  \mathsf{g}%
_{1}\right)  ,\Psi\left(  \mathsf{g}_{2}\right)  ,\ldots,\Psi\left(
\mathsf{g}_{n^{\prime}}\right)  \right]  ,\ \ \ \mathsf{g}_{j}\in
\mathsf{G}^{\left[  n^{\prime}\right]  },\ \Psi\left(  \mathsf{g}_{j}\right)
\in\mathsf{H}^{\left[  n^{\prime}\right]  },
\end{equation}
being algebra and grading group homomorphisms, respectively. In~addition, they
\textquotedblleft preserve\textquotedblright\ grading, such that%
\begin{equation}
\Phi\left(  \mathcal{A}\left(  \mathsf{g}\right)  \right)  \subseteq
\mathcal{B}\left(  \Psi\left(  \mathsf{g}\right)  \right)  ,\ \ \ \ \mathsf{g}%
\in\mathsf{G}^{\left[  n^{\prime}\right]  },\ \ \Psi\left(  \mathsf{g}\right)
\in\mathsf{H}^{\left[  n^{\prime}\right]  }, \label{fa}%
\end{equation}
for every graded component in (\ref{a}).
\end{Definition}

In the case $\Psi=\operatorname*{id}$, the~grading groups coincide
$\mathsf{G}^{\left[  n^{\prime}\right]  }=\mathsf{H}^{\left[  n^{\prime
}\right]  }$, and~homogeneous elements in $\mathcal{A}^{\left[  m,n\right]  }$
are really preserved as%
\begin{equation}
\Phi\left(  \mathcal{A}\left(  \mathsf{g}\right)  \right)  \subseteq
\mathcal{B}\left(  \mathsf{g}\right)  ,\ \ \ \ \mathsf{g}\in\mathsf{G}%
^{\left[  n^{\prime}\right]  },
\end{equation}
which means that they are mapped to the same corresponding homogeneous
elements   in $\mathcal{B}^{\left[  m,n\right]  }$.

\begin{Definition}
A graded homomorphism of polyadic algebras, which is defined by the pair
$\dbinom{\Phi}{\Psi}$ of maps $\Phi:\mathcal{A}^{\left[  m,n\right]
}\rightarrow\mathcal{B}^{\left[  m,n\right]  }$, $\Psi:\mathsf{G}^{\left[
n^{\prime}\right]  }\mathsf{\rightarrow H}^{\left[  n^{\prime}\right]  }$,
becomes a graded isomorphism (for binary case see, e.g.~\cite{bob/das/wyk}),
if both $\Phi$ and $\Psi$ are bijective, which preserves the grading by
definition of bijectivity, and~therefore, the inverse $\Phi^{-1}:\mathcal{B}%
^{\left[  m,n\right]  }\rightarrow\mathcal{A}^{\left[  m,n\right]  }$
$\Psi^{-1}:\mathsf{H}^{\left[  n^{\prime}\right]  }\mathsf{\rightarrow
G}^{\left[  n^{\prime}\right]  }$, is also a graded homomorphism $\dbinom
{\Phi^{-1}}{\Psi^{-1}}$. If~both inverse maps exist, the~polyadic graded
algebras $\mathcal{A}^{\left[  m,n\right]  }$ and $\mathcal{B}^{\left[
m,n\right]  }$ are called graded isomorphic. If~$\mathcal{A}^{\left[
m,n\right]  }=\mathcal{B}^{\left[  m,n\right]  }$, and~$\mathsf{G}^{\left[
n^{\prime}\right]  }=\mathsf{H}^{\left[  n^{\prime}\right]  }$, the~graded
isomorphism becomes a graded automorphism.
\end{Definition}

\begin{Example}
[Graded homomorphism]\label{ex-grad-hom}Let $\mathcal{A}^{\left[  2,3\right]
}$ be the ternary graded algebra from \textit{Example }\ref{ex-tern-pol}.
Consider the map $\Phi:$ $\mathcal{A}^{\left[  2,3\right]  }\rightarrow
\mathcal{A}^{\left[  2,3\right]  }$ defined by%
\begin{equation}
\Phi\left(  f\right)  \left(  x,y,z\right)  =f\left(  x+1,y+1,z+1\right)
,\ \ \ f\in\mathcal{A}^{\left[  2,3\right]  }. \label{f}%
\end{equation}

{The map $\Phi$ preserves (see (\ref{fm})) ternary multiplication $\Phi\left(
\mu^{\left[  3\right]  }\left[  f,g,h\right]  \right)  =\mu^{\left[  3\right]
}\left[  \Phi\left(  f\right)  ,\Phi\left(  g\right)  ,\Phi\left(  h\right)
\right]  $,} $f,g,h\in\mathcal{A}^{\left[  2,3\right]  }$, {since} %MDPI: Please check and confirm if the numbers within square brackets [2,3] refer to reference numbers. If yes, please write them as \cite{xyz} and move citation to text, instead of in equation.
%
%%%%%%Author's ANSWER: Do not move anything, because this is part of formula - arity. See explanation in Notation 1 and formula 3.
\begin{align}
&  \Phi\left(  \mu^{\left[  3\right]  }\left[  f,g,h\right]  \right)  \left(
x,y,z\right)  =\mu^{\left[  3\right]  }\left[  f,g,h\right]  \left(
x+1,y+1,z+1\right) \nonumber\\
&  =f\left(  x+1,y+1,z+1\right)  \cdot g\left(  x+1,y+1,z+1\right)  \cdot
h\left(  x+1,y+1,z+1\right)  ,
\end{align}
which is equal to%
\begin{align}
&  \mu^{\left[  3\right]  }\left[  \Phi\left(  f\right)  ,\Phi\left(
g\right)  ,\Phi\left(  h\right)  \right]  \left(  x,y,z\right) \nonumber\\
&  =f\left(  x+1,y+1,z+1\right)  \cdot g\left(  x+1,y+1,z+1\right)  \cdot
h\left(  x+1,y+1,z+1\right)  .
\end{align}

The map $\Phi$ (\ref{f}) preserves (see (\ref{fn})) the binary addition, which
follows from the ordinary additivity of polynomials $\Phi\left(  f+g\right)
=\Phi\left(  f\right)  +\Phi\left(  g\right)  $. Therefore, $\Phi$ is a
ternary homomorphism of the ternary algebra $\mathcal{A}^{\left[  2,3\right]
}$.

To show that $\Phi$ preserves grading, for~each graded component
$\mathcal{A}\left(  p\right)  $ in the decomposition (\ref{a3p}), we should
have $\Phi\left(  \mathcal{A}\left(  p\right)  \right)  \subseteq
\mathcal{A}\left(  p\right)  $ (\ref{fa}). Indeed, acting by $\Phi$ on any
monomial (\ref{pqr}) does not increase its degree, because~no terms of higher
degree appear, and~so, obviously, $\deg\left(  \Phi\left(  x^{p^{\prime}%
}y^{q^{\prime}}z^{r^{\prime}}\right)  \right)  =\deg\left(  \left(
x+1\right)  ^{p^{\prime}}\left(  y+1\right)  ^{q^{\prime}}\left(  z+1\right)
^{r^{\prime}}\right)  $, and~$\deg\Phi\left(  f\right)  =\deg f$, which means
that $\Phi$ maps each homogeneous component to~itself.

Thus, the~map $\Phi$ (\ref{f}) is a graded homomorphism (being of zero
grading) of the ternary algebra $\mathcal{A}^{\left[  2,3\right]  }$, it is
nontrivial since it changes polynomials; moreover, it is a graded
automorphism, because~$\Phi$ is an isomorphism of the algebra $\mathcal{A}%
^{\left[  2,3\right]  }\mathcal{\ }$to itself, and~$\Phi\circ\Phi^{-1}%
=\Phi^{-1}\circ\Phi=\operatorname*{id}$, where $\Phi^{-1}\left(  f\right)
\left(  x,y,z\right)  =f\left(  x-1,y-1,z-1\right)  $, $f\in\mathcal{A}%
^{\left[  2,3\right]  }$.
\end{Example}

The main difference of polyadic structures (rings, fields, algebras, groups)
from their binary counterparts is not the necessity of the existence of identities,
units (\ref{ea}) or zeroes (\ref{az}), such that there exist, e.g.,~unitless
and zeroless fields or algebras, as~well as polyadic groups without identity
or many identities. Instead, the~invertibility is governed by the existence of
querelements (\ref{aq}) for multiplications and additions (for review and
examples, see~\cite{duplij2022}). Therefore, the properties of graded
homomorphisms are changed correspondingly: mappings of unit to unit and zero
to zero are demanded, if~they exist. However, querelements should be
mapped to querelements. Thus, we~have

\begin{Proposition}
If and only if both algebras in $\Phi:\mathcal{A}^{\left[  m,n\right]  }\rightarrow
\mathcal{B}^{\left[  m,n\right]  }$ have units (\ref{ea}) $\mathbf{e}_{A}$ and
$\mathbf{e}_{B}$ (are unital), and~both polyadic groups in $\Psi
:\mathsf{G}^{\left[  n^{\prime}\right]  }\rightarrow\mathsf{H}^{\left[
n^{\prime}\right]  }$ have identities $\mathsf{e}_{G}$ and $\mathsf{e}_{H}$,
respectively, then%
\begin{align}
\Phi\left(  \mathbf{e}_{A}\right)   &  =\mathbf{e}_{B},\\
\Psi\left(  \mathsf{e}_{G}\right)   &  =\mathsf{e}_{H}.
\end{align}

\end{Proposition}

\begin{proof}
This follows directly from application of the mapping $\Phi,\Psi$ to the
definitions (\ref{ea}) and using the homomorphism property.
\end{proof}

Obviously,

\begin{Assertion}
Sets of elements in $\mathcal{A}^{\left[  m,n\right]  }$ and $\mathsf{G}%
^{\left[  n^{\prime}\right]  }$ that map under $\Phi$ and $\Psi$ to the unit
and identity of $\mathcal{B}^{\left[  m,n\right]  }$ and $\mathsf{H}^{\left[
n^{\prime}\right]  }$ form subalgebra of $\mathcal{A}^{\left[  m,n\right]  }$
and subgroup of $\mathsf{G}^{\left[  n^{\prime}\right]  }$, respectively.
\end{Assertion}

\begin{Proposition}
If and only if both algebras in $\Phi:\mathcal{A}^{\left[  m,n\right]  }\rightarrow
\mathcal{B}^{\left[  m,n\right]  }$ have zeroes (\ref{ea}) $\mathbf{z}_{A}$
and $\mathbf{z}_{B}$,   respectively, then%
\begin{equation}
\Phi\left(  \mathbf{z}_{A}\right)  =\mathbf{z}_{B}.\label{fzz}%
\end{equation}

\end{Proposition}

\begin{proof}
This follows from (\ref{az}) and (\ref{fm}) directly.
\end{proof}

Denote $\mathbf{\bar{a}}$, $\mathbf{\bar{b}}$ and $\mathsf{\bar{g}}$,
$\mathsf{\bar{h}}$ as the querelements of polyadic algebra elements
$\mathbf{a}\in\mathcal{A}^{\left[  m,n\right]  }$, $\mathbf{b\in}%
\mathcal{B}^{\left[  m,n\right]  }$ and the grading $n^{\prime}$-ary group
elements $\mathsf{g}\in\mathsf{G}^{\left[  n^{\prime}\right]  }$,
$\mathsf{h}\in\mathsf{H}^{\left[  n^{\prime}\right]  }$, respectively   (see
(\ref{aq})).~Then,

\begin{Proposition}
The querelements are mapped to the corresponding querelements, as%
\begin{align}
\Phi\left(  \mathbf{a}\right)   &  =\mathbf{b}\Longrightarrow\Phi\left(
\mathbf{\bar{a}}\right)  =\mathbf{\bar{b}},\ \ \ \mathbf{a}\in\mathcal{A}%
^{\left[  m,n\right]  },\ \mathbf{b\in}\mathcal{B}^{\left[  m,n\right]  }\\
\Psi\left(  \mathsf{g}\right)   &  =\mathsf{g}\Longrightarrow\Psi\left(
\mathsf{\bar{g}}\right)  =\mathsf{\bar{h}},\ \ \ \mathsf{g}\in\mathsf{G}%
^{\left[  n^{\prime}\right]  },\ \mathsf{h}\in\mathsf{H}^{\left[  n^{\prime
}\right]  }.
\end{align}

\end{Proposition}

\begin{proof}
This follows directly from (\ref{aq}) and (\ref{fm}).
\end{proof}

\begin{Corollary}
If $\mathcal{A}_{sub}^{\left[  m,n\right]  }$ is a graded subalgebra of
$\mathcal{A}^{\left[  m,n\right]  }$, then $\Phi\left(  \mathcal{A}%
_{sub}^{\left[  m,n\right]  }\right)  $ is a graded subalgebra of $\Phi\left(
\mathcal{A}^{\left[  m,n\right]  }\right)  $. Oppositely, if~$\mathcal{B}%
_{sub}^{\left[  m,n\right]  }$ is a graded subalgebra of $\mathcal{B}^{\left[
m,n\right]  }$, then $\Phi^{-1}\left(  \mathcal{B}_{sub}^{\left[  m,n\right]
}\right)  $ is a graded subalgebra of $\Phi^{-1}\left(  \mathcal{B}^{\left[
m,n\right]  }\right)  $.
\end{Corollary}

\begin{Definition}
If $\Phi:\mathcal{A}^{\left[  m,n\right]  }\rightarrow\mathcal{B}^{\left[
m,n\right]  }$ is the graded homomorphism, then its image is defined by%
\begin{equation}
\operatorname*{im}\Phi=\left\{  \Phi\left(  \mathbf{a}\right)  \mid
\mathbf{a}\in\mathcal{A}^{\left[  m,n\right]  }\right\}  \subseteq
\mathcal{B}^{\left[  m,n\right]  }.
\end{equation}

\end{Definition}

Because $\mathcal{A}^{\left[  m,n\right]  }$ and $\mathcal{B}^{\left[
m,n\right]  }$ are graded having the direct sum decompositions (\ref{a})   and
(\ref{b}), the~subspaces $\operatorname*{im}\Phi\cap\mathcal{B}\left(
\mathsf{g}_{i}\right)  =\Phi\left(  \mathcal{A}\left(  \mathsf{g}_{i}\right)
\right)  $ form a grading (decomposition) of the image, as~follows:%
\begin{equation}
\operatorname*{im}\Phi=\bigoplus_{\mathsf{g}_{i}\in\mathsf{G}^{\left[
n^{\prime}\right]  }}\Phi\left(  \mathcal{A}\left(  \mathsf{g}_{i}\right)
\right)  ,
\end{equation}
which is a polyadic graded subalgebra, since it follows from (\ref{ma}) and the
homomorphism multiplicative property of $\Phi$, that%
\begin{equation}
\mathbf{\mu}_{b}^{\left[  n\right]  }\left[  \Phi\left(  \mathcal{A}\left(
\mathsf{g}_{1}\right)  \right)  ,\Phi\left(  \mathcal{A}\left(  \mathsf{g}%
_{2}\right)  \right)  ,\ldots,\Phi\left(  \mathcal{A}\left(  \mathsf{g}%
_{n}\right)  \right)  \right]  \subseteq\Phi\left(  \mathcal{A}\left(  \mu
_{g}^{\left[  n^{\prime}\right]  }\left[  \mathsf{g}_{1},\mathsf{g}_{2}%
,\ldots,\mathsf{g}_{n^{\prime}}\right]  \right)  \right)  . \label{mf}%
\end{equation}

\begin{Remark}
\label{rem-ker}In contrast to the binary case, where a zero element always
exists, and the kernel is therefore always definable as the preimage of zero,
in the polyadic setting, the standard kernel can be defined only if the
polyadic algebra $\mathcal{B}^{\left[  m,n\right]  }$ possesses a zero.
\end{Remark}

\begin{Definition}
[Standard kernel]If and only if the algebra $\mathcal{B}^{\left[  m,n\right]
}$ in $\Phi:\mathcal{A}^{\left[  m,n\right]  }\rightarrow\mathcal{B}^{\left[
m,n\right]  }$ has zero (\ref{ea}) $\mathbf{z}_{B}$, then the standard kernel
of $\Phi$ can be defined by%
\begin{equation}
\ker\Phi=\left\{  \mathbf{a}\in\mathcal{A}^{\left[  m,n\right]  }\mid
\Phi\left(  \mathbf{a}\right)  =\mathbf{z}_{B}\right\}  . \label{kf}%
\end{equation}

\end{Definition}

Because $\Phi$ is a polyadic algebra homomorphism (\ref{fm}), its standard
kernel (\ref{kf}) is an ideal of $\mathcal{A}^{\left[  m,n\right]  }$ (not
only a graded subalgebra, as~$\operatorname*{im}\Phi$). Moreover, $\ker\Phi$
inherits a grading from $\mathcal{A}^{\left[  m,n\right]  }$ becoming a graded
ideal, which means that there is decomposition on homogeneous components%
\begin{equation}
\ker\Phi=\bigoplus_{\mathsf{g}_{i}\in\mathsf{G}^{\left[  n^{\prime}\right]  }%
}\left(  \ker\Phi\cap\mathcal{A}\left(  \mathsf{g}_{i}\right)  \right)  ,
\label{ka}%
\end{equation}
such that, if~an element belongs to the standard kernel, than~each of its
homogeneous components also belongs to it. Thus, $\ker\Phi$ is always
generated, as~an ideal, by~the homogeneous components of any element of the
standard kernel (\ref{kf}).

A more general definition of a kernel is needed to avoid difficulty with the absence
of some elements; see \textit{Remark} \ref{rem-ker}. In~universal algebra~\cite{cohn,lang} the kernel is defined as a special equivalence relation
(congruence) on the domain and~not as a set of elements mapping to the
identity (for groups) or to zero for $\mathcal{B}^{\left[  m,n\right]  }$
(\ref{kf}).

\begin{Definition}
[Congruence kernel]The congruence kernel of the homomorphism $\Phi$, denoted
$\ker_{\theta}\Phi$, is the equivalence relation%
\begin{equation}
\theta\equiv\ker_{\theta}\Phi=\left\{  \left(  \mathbf{a}_{1},\mathbf{a}%
_{2}\right)  \in\mathcal{A}^{\left[  m,n\right]  }\times\mathcal{A}^{\left[
m,n\right]  }\mid\Phi\left(  \mathbf{a}_{1}\right)  =\Phi\left(
\mathbf{a}_{2}\right)  \right\}  , \label{kc}%
\end{equation}
which is compatible with all operations in $\mathcal{A}^{\left[  m,n\right]
}$.
\end{Definition}

The set definition (\ref{kf}) is a special case of the congruence definition
(\ref{kc}), where the algebraic structure gives the possibility (which is not
always the case, see {Remark} \ref{rem-ker}) to describe the congruence
by a single distinguished congruence class: the preimage of identity or zero.
The congruence $\ker_{\theta}\Phi$ respects grading, such that it only relates
elements of the same degree (cf. (\ref{ka})):%
\begin{equation}
\ker_{\theta}\Phi=\bigoplus_{\mathsf{g}_{i}\in\mathsf{G}^{\left[  n^{\prime
}\right]  }}\left(  \ker_{\theta}\Phi\cap\mathcal{A}\left(  \mathsf{g}%
_{i}\right)  \times\mathcal{A}\left(  \mathsf{g}_{i}\right)  \right)  .
\end{equation}

The congruence definition (\ref{kc}) allows us to form the quotient polyadic
algebra $\mathcal{A}^{\left[  m,n\right]  }\diagup\ker_{\theta}\Phi$, whose
elements are the equivalence classes of the congruence $\theta$ denoted by
$\lfloor\mathbf{a\rfloor}_{\theta}$ for $\mathbf{a}\in\mathcal{A}^{\left[
m,n\right]  }$:%
\begin{equation}
\mathcal{A}^{\left[  m,n\right]  }\diagup\ker_{\theta}\Phi=\mathcal{\bar{A}%
}^{\left[  m,n\right]  }=\left\langle \left\{  \lfloor\mathbf{a\rfloor
}_{\theta}\right\}  \mid\mathbf{\bar{\nu}}^{\left[  m\right]  },\mathbf{\bar
{\mu}}^{\left[  n\right]  }\right\rangle ,\ \ \ \mathbf{a}\in\mathcal{A}%
^{\left[  m,n\right]  }. \label{ak}%
\end{equation}

The operations on classes $\mathbf{\bar{\nu}}^{\left[  m\right]
},\mathbf{\bar{\mu}}^{\left[  n\right]  }$ can be naturally defined using
representatives,   as follows:%
\begin{align}
\mathbf{\bar{\nu}}_{a}^{\left[  m\right]  }\left[  \lfloor\mathbf{a}%
_{1}\mathbf{\rfloor}_{\theta},\lfloor\mathbf{a}_{2}\mathbf{\rfloor}_{\theta
},\ldots,\lfloor\mathbf{a}_{m}\mathbf{\rfloor}_{\theta}\right]   &
\equiv\lfloor\mathbf{\nu}_{a}^{\left[  m\right]  }\left[  \mathbf{a}%
_{1},\mathbf{a}_{2},\ldots,\mathbf{a}_{m}\right]  \rfloor_{\theta},\\
\mathbf{\bar{\mu}}^{\left[  n\right]  }\left[  \lfloor\mathbf{a}%
_{1}\mathbf{\rfloor}_{\theta},\lfloor\mathbf{a}_{2}\mathbf{\rfloor}_{\theta
},\ldots,\lfloor\mathbf{a}_{m}\mathbf{\rfloor}_{\theta}\right]   &
\equiv\lfloor\mathbf{\mu}_{a}^{\left[  n\right]  }\left[  \mathbf{a}%
_{1},\mathbf{a}_{2},\ldots,\mathbf{a}_{n}\right]  \rfloor_{\theta
},\ \ \ \mathbf{a}_{i}\in\mathcal{A}^{\left[  m,n\right]  },
\end{align}
and they are well-defined, because~the resulting r.h.s. do not depend on the
choice of the representatives (which directly follows from the homomorphism
property and (\ref{kc})). The~quotient polyadic algebra $\mathcal{A}^{\left[
m,n\right]  }\diagup\ker_{\theta}\Phi$ respects multiary grading, such that
the decomposition (\ref{a}) with the compatibility (\ref{ma}) provides
\begin{equation}
\mathcal{\bar{A}}^{\left[  m,n\right]  }=\bigoplus_{\mathsf{g}_{i}%
\in\mathsf{G}^{\left[  n^{\prime}\right]  }}\mathcal{\bar{A}}\left(
\mathsf{g}_{i}\right)  ,
\end{equation}
and therefore, $\mathcal{\bar{A}}^{\left[  m,n\right]  }$ is a multiary graded
polyadic~algebra.

\begin{Theorem}
[First Isomorphism Theorem]Let $\mathcal{A}^{\left[  m,n\right]  }$ and
$\mathcal{B}^{\left[  m,n\right]  }$ be multiary graded polyadic algebras, and~
let $\Phi:\mathcal{A}^{\left[  m,n\right]  }\rightarrow\mathcal{B}^{\left[
m,n\right]  }$ be a multiary graded homomorphism; then, the quotient algebra
(\ref{ak}) is isomorphic to the graded image $\operatorname*{im}\Phi$ being a
subalgebra of $\mathcal{B}^{\left[  m,n\right]  }$,
\begin{equation}
\mathcal{A}^{\left[  m,n\right]  }\diagup\ker_{\theta}\Phi\cong%
\operatorname*{im}\Phi. \label{akf}%
\end{equation}

\end{Theorem}

\begin{proof}
Define the map $\bar{\Phi}$ on the congruence classes (\ref{ak})
$\mathcal{A}^{\left[  m,n\right]  }\diagup\ker_{\theta}\Phi\rightarrow
\Phi\left(  \mathcal{A}^{\left[  m,n\right]  }\right)  $ by%
\begin{equation}
\bar{\Phi}\left(  \lfloor\mathbf{a\rfloor}_{\theta}\right)  =\Phi\left(
\mathbf{a}\right)  ,\ \ \ \mathbf{a}\in\mathcal{A}^{\left[  m,n\right]  }.
\label{ff}%
\end{equation}

It is well-defined, because, if~$\lfloor\mathbf{a\rfloor}_{\theta}%
=\lfloor\mathbf{b\rfloor}_{\theta}$, $\mathbf{a},\mathbf{b}\in\mathcal{A}%
^{\left[  m,n\right]  }$, then $\left(  \mathbf{a},\mathbf{b}\right)  \in
\ker_{\theta}\Phi$, and~so $\Phi\left(  \mathbf{a}\right)  =\Phi\left(
\mathbf{b}\right)  $; thus, $\bar{\Phi}\left(  \lfloor\mathbf{a\rfloor}%
_{\theta}\right)  =\bar{\Phi}\left(  \lfloor\mathbf{b\rfloor}_{\theta}\right)
$.

The map $\bar{\Phi}$ respects $\mathsf{G}$-grading: if $\lfloor
\mathbf{a\rfloor}_{\theta}\in\mathcal{\bar{A}}\left(  \mathsf{g}\right)  $,
$\mathsf{g}\in\mathsf{G}^{\left[  n^{\prime}\right]  }$, then $\mathbf{a}%
\in\mathcal{A}\left(  \mathsf{g}\right)  $, and~so $\Phi\left(  \mathbf{a}%
\right)  \in\mathcal{B}\left(  \mathsf{g}\right)  \cap\Phi\left(
\mathcal{A}^{\left[  m,n\right]  }\right)  =\Phi\left(  \mathcal{A}\left(
\mathsf{g}\right)  \right)  $, the~graded component of the image
$\operatorname*{im}\Phi$.

It is an algebra map that preserves the polyadic operations; that is, we derive%
%\vspace{-15pt}
%\begin{adjustwidth}{-\extralength}{0cm}
\begin{align}
&  \bar{\Phi}\left(  \mathbf{\bar{\nu}}_{a}^{\left[  m\right]  }\left[
\lfloor\mathbf{a}_{1}\mathbf{\rfloor}_{\theta},\lfloor\mathbf{a}%
_{2}\mathbf{\rfloor}_{\theta},\ldots,\lfloor\mathbf{a}_{m}\mathbf{\rfloor
}_{\theta}\right]  \right)  =\bar{\Phi}\left(  \lfloor\mathbf{\nu}%
_{a}^{\left[  m\right]  }\left[  \mathbf{a}_{1},\mathbf{a}_{2},\ldots
,\mathbf{a}_{m}\right]  \rfloor_{\theta}\right)  =\Phi\left(  \mathbf{\nu}%
_{a}^{\left[  m\right]  }\left[  \mathbf{a}_{1},\mathbf{a}_{2},\ldots
,\mathbf{a}_{m}\right]  \right) \\
&  =\mathbf{\nu}_{a}^{\left[  m\right]  }\left[  \Phi\left(  \mathbf{a}%
_{1}\right)  ,\Phi\left(  \mathbf{a}_{2}\right)  ,\ldots,\Phi\left(
\mathbf{a}_{m}\right)  \right]  =\mathbf{\nu}_{a}^{\left[  m\right]  }\left[
\bar{\Phi}\left(  \lfloor\mathbf{a}_{1}\mathbf{\rfloor}_{\theta}\right)
,\bar{\Phi}\left(  \lfloor\mathbf{a}_{2}\mathbf{\rfloor}_{\theta}\right)
,\ldots,\bar{\Phi}\left(  \lfloor\mathbf{a}_{m}\mathbf{\rfloor}_{\theta
}\right)  \right]  ,
\end{align}%
%\end{adjustwidth}
%\vspace{-25pt}
%\begin{adjustwidth}{-\extralength}{0cm}
\begin{align}
&  \bar{\Phi}\left(  \mathbf{\bar{\mu}}^{\left[  n\right]  }\left[
\lfloor\mathbf{a}_{1}\mathbf{\rfloor}_{\theta},\lfloor\mathbf{a}%
_{2}\mathbf{\rfloor}_{\theta},\ldots,\lfloor\mathbf{a}_{m}\mathbf{\rfloor
}_{\theta}\right]  \right)  =\bar{\Phi}\left(  \lfloor\mathbf{\mu}%
_{a}^{\left[  n\right]  }\left[  \mathbf{a}_{1},\mathbf{a}_{2},\ldots
,\mathbf{a}_{n}\right]  \rfloor_{\theta}\right)  =\Phi\left(  \mathbf{\mu}%
_{a}^{\left[  n\right]  }\left[  \mathbf{a}_{1},\mathbf{a}_{2},\ldots
,\mathbf{a}_{n}\right]  \right) \\
&  =\mathbf{\mu}_{a}^{\left[  n\right]  }\left[  \Phi\left(  \mathbf{a}%
_{1}\right)  ,\Phi\left(  \mathbf{a}_{2}\right)  ,\ldots,\Phi\left(
\mathbf{a}_{m}\right)  \right]  =\mathbf{\mu}_{a}^{\left[  n\right]  }\left[
\bar{\Phi}\left(  \lfloor\mathbf{a}_{1}\mathbf{\rfloor}_{\theta}\right)
,\bar{\Phi}\left(  \lfloor\mathbf{a}_{2}\mathbf{\rfloor}_{\theta}\right)
,\ldots,\bar{\Phi}\left(  \lfloor\mathbf{a}_{n}\mathbf{\rfloor}_{\theta
}\right)  \right]  ,\ \ \ \mathbf{a}_{i}\in\mathcal{A}^{\left[  m,n\right]  },
\end{align}
%\end{adjustwidth}
because $\Phi$ is a homomorphism of $\mathcal{A}^{\left[  m,n\right]  }$.

The map $\bar{\Phi}$ is injective, such that, if~$\bar{\Phi}\left(
\lfloor\mathbf{a\rfloor}_{\theta}\right)  =\bar{\Phi}\left(  \lfloor
\mathbf{b\rfloor}_{\theta}\right)  $, $\mathbf{a},\mathbf{b}\in\mathcal{A}%
^{\left[  m,n\right]  }$, then $\Phi\left(  \mathbf{a}\right)  =\Phi\left(
\mathbf{b}\right)  $, and~so, $\left(  \mathbf{a},\mathbf{b}\right)  \in
\ker_{\theta}\Phi$; hence, $\lfloor\mathbf{a\rfloor}_{\theta}=\lfloor
\mathbf{b\rfloor}_{\theta}$.

It is surjective, because~for any $\mathbf{b}\in\Phi\left(  A\right)  $, there
exists $\mathbf{a}\in A$ with $\mathbf{b}=\Phi\left(  \mathbf{a}\right)  $,
and so, $\bar{\Phi}\left(  \lfloor\mathbf{a\rfloor}_{\theta}\right)
=\mathbf{b}$.

Therefore, the~map $\bar{\Phi}$ (\ref{ff}) is a surjective and injective (i.e.
bijective) $\mathsf{G}$-graded homomorphism, hence an isomorphism (\ref{akf}).
\end{proof}

Thus, the~First Isomorphism Theorem for $\mathsf{G}$-graded polyadic algebras
states that the structure of the image of a homomorphism is exactly the same
as the structure of the domain ``collapsed'' by the equivalence relation induced
by the map. This formulation ensures that all structural properties (the
polyadic algebra operations and the multiary group grading) are rigorously
preserved through the~isomorphism.

Let us denote the set of $N$ elements by $G_{N}=\left\{  0,1,\ldots,N\right\}
\in\mathbb{N}_{0}$. The~simplest binary group on $G_{2}=\left\{
\mathsf{0,1}\right\}  $ is the cyclic group of order $2$ (or integers modulo
$2$) $\mathbb{Z}_{2}=\mathsf{G}^{\left[  2^{\prime}\right]  }=\left\langle
G_{2}\mid\mu_{g}^{\left[  2^{\prime}\right]  }\right\rangle $ with the binary
operation $\mu_{g}^{\left[  2^{\prime}\right]  }\left[  \mathsf{x}%
,\mathsf{y}\right]  =\left(  \mathsf{x}+\mathsf{y}\right)  \operatorname{mod}%
2$, $\mathsf{x},\mathsf{y}\in G_{2}$.

\begin{Example}
Consider the polynomial ring $\mathbb{Z}\left[  t\right]  $ graded by degree,
reduced $\operatorname{mod}6$ and $\operatorname{mod}3$. The~corresponding
$\mathbb{Z}$-graded ternary algebra $\mathcal{A}^{\left[  3,3\right]  }$ with
zero and the zeroless $\mathbb{Z}$-graded ternary algebra $\mathcal{B}%
^{\left[  3,3\right]  }$ are%

%\vspace{-15pt}
%\begin{adjustwidth}{-\extralength}{0cm}
\begin{align}
\mathcal{A}^{\left[  3,3\right]  }  &  =\left\langle \mathbb{Z}\left[
t\right]  \mid\mathbf{\nu}_{a}^{\left[  3\right]  },\mathbf{\mu}_{a}^{\left[
3\right]  },\mathbf{\mu}_{Z}^{\left[  0\right]  }\right\rangle =\bigoplus
_{k\geq0}\mathcal{A}\left(  k\right)  ,\mathcal{A}\left(  k\right)  =\left\{
a_{k}\cdot t^{k}\right\}  ,a_{k}\in\mathbb{Z\diagup}6\mathbb{Z}=0^{\prime
},1^{\prime},2^{\prime},3^{\prime},4^{\prime},5^{\prime},\label{a33}\\
\mathcal{B}^{\left[  3,3\right]  }  &  =\left\langle \mathbb{Z}\left[
t\right]  \mid\mathbf{\nu}_{b}^{\left[  3\right]  },\mathbf{\mu}_{b}^{\left[
3\right]  }\right\rangle =\bigoplus_{k\geq0}\mathcal{B}\left(  k\right)
,\ \ \ \mathcal{B}\left(  k\right)  =\left\{  b_{k}\cdot t^{k}\right\}
,\ \ \ b_{k}\in\mathbb{Z\diagup}3\mathbb{Z}=0^{\prime\prime},1^{\prime\prime
},2^{\prime\prime}, \label{b33}%
\end{align}
%\end{adjustwidth}
where $\mathcal{A}^{\left[  3,3\right]  }$ is a ternary algebra with zero
named in its signature as the null-ary operation $\mathbf{\mu}_{Z}^{\left[
0\right]  }=\mathbf{z}_{A}=0\in\mathbb{Z\diagup}6\mathbb{Z}$, and
$\mathcal{B}^{\left[  3,3\right]  }$ is a ternary algebra without zero (in
signature), without binary substraction. The~coefficients $a_{k}$ and
$b_{k}$ are in (different) congruence classes of $\mathbb{Z\diagup}%
6\mathbb{Z}$ and $\mathbb{Z\diagup}3\mathbb{Z}$, respectively. The~graded
components $\mathcal{A}\left(  k\right)  $ ($\mathcal{B}\left(  k\right)  $)
have $\deg(\mathcal{A}\left(  k\right)  )=\deg(\mathcal{B}\left(  k\right)
)=k$. The~operations in $\mathcal{A}^{\left[  3,3\right]  }$ and
$\mathcal{B}^{\left[  3,3\right]  }$ are%
\begin{align}
\mathbf{\nu}_{a}^{\left[  3\right]  }\left[  x,y,z\right]   &  =\left(
x+y+z\right)  \operatorname{mod}6\in\mathbb{Z\diagup}6\mathbb{Z}%
,\ \ \ \mathbf{\mu}_{a}^{\left[  3\right]  }\left[  x,y,z\right]  =\left(
xyz\right)  \operatorname{mod}6\in\mathbb{Z\diagup}6\mathbb{Z},\nonumber\\
\mathbf{\mu}_{a,Z}^{\left[  0\right]  }  &  =0\in\mathbb{Z\diagup}%
6\mathbb{Z},\label{nu3}\\
\mathbf{\nu}_{b}^{\left[  3\right]  }\left[  x,y,z\right]   &  =\left(
x+y+z\right)  \operatorname{mod}3\in\mathbb{Z\diagup}3\mathbb{Z}%
,\ \ \ \mathbf{\mu}_{b}^{\left[  3\right]  }\left[  x,y,z\right]  =\left(
xyz\right)  \operatorname{mod}3\in\mathbb{Z\diagup}3\mathbb{Z}. \label{mu3}%
\end{align}

By construction, the~ternary operations (\ref{nu3}) and (\ref{mu3}) satisfy
ternary distributivity and are totally associative and commutative, as
both algebras are, and they respect grading by the group $\mathsf{G}%
=\mathbb{Z}_{\geq0}$; that is,%
%\vspace{-12pt}
%\begin{adjustwidth}{-\extralength}{0cm}
\begin{align}
\mathbf{\nu}_{a}^{\left[  3\right]  }\left[  \mathcal{A}\left(  k_{1}\right)
\mathcal{A}\left(  k_{2}\right)  \mathcal{A}\left(  k_{3}\right)  \right]   &
\subseteq\mathcal{A}\left(  k_{1}+k_{2}+k_{3}\right)  ,\ \ \ \mathbf{\mu}%
_{a}^{\left[  3\right]  }\left[  \mathcal{A}\left(  k_{1}\right)
\mathcal{A}\left(  k_{2}\right)  \mathcal{A}\left(  k_{3}\right)  \right]
\subseteq\mathcal{A}\left(  k_{1}+k_{2}+k_{3}\right)  ,\\
\mathbf{\nu}_{b}^{\left[  3\right]  }\left[  \mathcal{B}\left(  k_{1}\right)
\mathcal{B}\left(  k_{2}\right)  \mathcal{B}\left(  k_{3}\right)  \right]   &
\subseteq\mathcal{B}\left(  k_{1}+k_{2}+k_{3}\right)  ,\ \ \ \mathbf{\mu}%
_{b}^{\left[  3\right]  }\left[  \mathcal{B}\left(  k_{1}\right)
\mathcal{B}\left(  k_{2}\right)  \mathcal{B}\left(  k_{3}\right)  \right]
\subseteq\mathcal{B}\left(  k_{1}+k_{2}+k_{3}\right)  .
\end{align}
%\end{adjustwidth}

We define the homomorphism $\Phi:\mathcal{A}^{\left[  3,3\right]  }%
\rightarrow\mathcal{B}^{\left[  3,3\right]  }$ by%
\begin{equation}
\Phi\left(  x\right)  =x\operatorname{mod}3.\label{fx}%
\end{equation}

The map $\Phi$ is a graded homomorphism (of degree zero), because~it acts
degreewise $\Phi\left(  a\cdot t^{k}\right)  =\left(  a\operatorname{mod}%
3\right)  \cdot t^{k}$, and~therefore, $\Phi\left(  \mathcal{A}\left(
k\right)  \right)  \subseteq\mathcal{B}\left(  k\right)  $. The~standard
kernel (\ref{kf}) cannot be defined, because~the target algebra $\mathcal{B}%
^{\left[  3,3\right]  }$ is zeroless (no zero as null-operation $\mathbf{z}%
_{B}$ in its signature (\ref{b33}), and~there is no relation (\ref{fzz}) for
$\Phi$). However, the~graded congruence kernel (\ref{kc}) can be defined as
$\theta=\ker_{\theta}\Phi$ by the equivalence relation%
\begin{equation}
x\theta y\Longleftrightarrow\Phi\left(  x\right)  =\Phi\left(  y\right)
\text{ and }\deg x=\deg y,\ \ \ x,y\in\mathbb{Z\diagup}6\mathbb{Z}.
\end{equation}

Using (\ref{kc}) and (\ref{fx}), we obtain three congruence classes at each
degree $k$%

%\begin{adjustwidth}{-\extralength}{0cm}
\begin{equation}
\lfloor0^{\prime}\cdot t^{k}\mathbf{\rfloor}_{\theta}=\left\{  0^{\prime
\prime}\cdot t^{k},3^{\prime\prime}\cdot t^{k}\right\}  ,\ \ \lfloor1^{\prime
}\cdot t^{k}\mathbf{\rfloor}_{\theta}=\left\{  1^{\prime\prime}\cdot
t^{k},4^{\prime\prime}\cdot t^{k}\right\}  ,\ \ \lfloor2^{\prime}\cdot
t^{k}\mathbf{\rfloor}_{\theta}=\left\{  2^{\prime\prime}\cdot t^{k}%
,5^{\prime\prime}\cdot t^{k}\right\}  ,
\end{equation}
%\end{adjustwidth}
and so the graded quotient becomes%
\begin{equation}
\mathcal{A}^{\left[  3,3\right]  }\diagup\ker_{\theta}\Phi=\bigoplus_{k\geq
0}\mathcal{A}\left(  k\right)  \diagup\theta,\ \ \ \ \ \ \ \mathcal{A}\left(
k\right)  \diagup\theta\cong\mathbb{Z\diagup}3\mathbb{Z}.
\end{equation}

Therefore, there is the isomorphism of $\mathbb{Z}$-graded ternary algebras
\begin{equation}
\mathcal{A}^{\left[  3,3\right]  }\diagup\ker_{\theta}\Phi\cong%
\operatorname*{im}\Phi,
\end{equation}
which preserves both the ternary algebra structure and the grading simultaneously.
\end{Example}

\section{\label{sec-super}\textsc{Ternary~Superalgebras}}

Let us consider examples of multiary ($n^{\prime}$-ary) abelian grading groups
$\mathsf{G}^{\left[  n^{\prime}\right]  }$ of lowest order $N$ in additive~notation.

\begin{Example}
[Derived ternary superalgebra]The simplest example is the derived ternary
algebra over the field $\Bbbk=\mathbb{R}$ on two elements, which is derived
from the binary algebra  $\mathcal{A}^{\left[  2,2\right]  }=\left\langle
\left\{  \mathbf{0},\mathbf{1}\right\}  \mid\mathbf{\nu}_{a}^{\left[
2\right]  }\equiv\left(  +\right)  ,\mathbf{\mu}_{a}^{\left[  2\right]
}\equiv\left(  \cdot\right)  \right\rangle $, as~follows:%
\begin{equation}
\mathcal{A}^{\left[  2,3\right]  }=\left\langle \left\{  \mathbf{0}%
,\mathbf{1}\right\}  \mid\mathbf{\nu}_{a}^{\left[  2\right]  }\equiv\left(
+\right)  ,\mathbf{\mu}_{a}^{\left[  3\right]  }\equiv\mathbf{\mu}%
_{a}^{\left[  2\right]  \circ2}\right\rangle,
\end{equation}
having the relations%
\begin{equation}
\mathbf{\mu}_{a}^{\left[  3\right]  }\left[  \mathbf{0},\mathbf{0}%
,\mathbf{0}\right]  =\mathbf{0},\ \ \ \mathbf{\mu}_{a}^{\left[  3\right]
}\left[  \mathbf{0},\mathbf{0},\mathbf{1}\right]  =\mathbf{1}%
,\ \ \ \mathbf{\mu}_{a}^{\left[  3\right]  }\left[  \mathbf{0},\mathbf{1}%
,\mathbf{1}\right]  =\mathbf{0},\ \ \ \mathbf{\mu}_{a}^{\left[  3\right]
}\left[  \mathbf{1},\mathbf{1},\mathbf{1}\right]  =\mathbf{0}. \label{m0}%
\end{equation}

The direct sum decomposition (\ref{a}) now is%
\begin{equation}
\mathcal{A}^{\left[  2,3\right]  }=\mathcal{A}\left(  \mathsf{0}\right)
\oplus\mathcal{A}\left(  \mathsf{1}\right)  ,
\end{equation}
where $\mathcal{A}\left(  \mathsf{0}\right)  =a\cdot\mathsf{0}$,
$\mathcal{A}\left(  \mathsf{1}\right)  =b\cdot\mathsf{1}$, $a,b\in\Bbbk$, and~
$\mathsf{0,1}\in\mathsf{G}^{\left[  2^{\prime}\right]  }\equiv\mathbb{Z}_{2}$
is the binary grading~group.

Therefore, $\mathcal{A}^{\left[  2,3\right]  }$ is $\mathbb{Z}_{2}$-graded
derived ternary algebra or~commutative ternary superalgebra:%
\begin{equation}%
\begin{array}
[c]{cc}%
\mathbf{\mu}_{a}^{\left[  3\right]  }\left[  \mathcal{A}\left(  0\right)
,\mathcal{A}\left(  0\right)  ,\mathcal{A}\left(  0\right)  \right]
\subseteq\mathcal{A}\left(  0\right)  ,\ \ \ \ \  & \mathbf{\mu}_{a}^{\left[
3\right]  }\left[  \mathcal{A}\left(  0\right)  ,\mathcal{A}\left(  0\right)
,\mathcal{A}\left(  1\right)  \right]  \subseteq\mathcal{A}\left(  1\right)
,\\
\mathbf{\mu}_{a}^{\left[  3\right]  }\left[  \mathcal{A}\left(  0\right)
,\mathcal{A}\left(  1\right)  ,\mathcal{A}\left(  1\right)  \right]
\subseteq\mathcal{A}\left(  0\right)  ,\ \ \ \ \  & \mathbf{\mu}_{a}^{\left[
3\right]  }\left[  \mathcal{A}\left(  1\right)  ,\mathcal{A}\left(  1\right)
,\mathcal{A}\left(  1\right)  \right]  \subseteq\mathcal{A}\left(  0\right)  .
\end{array}
\label{m3}%
\end{equation}

Define the mappings%
\begin{align}
\Phi\left(  \mathbf{0}\right)   &  =\mathbf{0},\ \ \ \ \Phi\left(
\mathbf{1}\right)  =c\cdot\mathbf{1},\ \ \ \ \ \ c\in\Bbbk,\ \ \ \ \mathbf{0}%
,\mathbf{1}\in\mathcal{A}^{\left[  2,3\right]  },\label{f0}\\
\Psi &  =\operatorname*{id}. \label{f1}%
\end{align}

The mapping $\Phi$ preserves the grading (see (\ref{fa}))%
\begin{equation}
\Phi\left(  \mathcal{A}\left(  \mathsf{0}\right)  \right)  \subseteq
\mathcal{A}\left(  \mathsf{0}\right)  ,\ \ \ \Phi\left(  \mathcal{A}\left(
\mathsf{1}\right)  \right)  \subseteq\mathcal{A}\left(  \mathsf{1}\right)  .
\label{fa0}%
\end{equation}

The only one nontrivial relation to show this follows from the second one in
(\ref{m0}):%
\begin{align}
\Phi\left(  \mathbf{\mu}_{a}^{\left[  3\right]  }\left[  \mathbf{0}%
,\mathbf{0},\mathbf{1}\right]  \right)   &  =\Phi\left(  \mathbf{1}\right)
=c\cdot\mathbf{1},\\
\mathbf{\mu}_{a}^{\left[  3\right]  }\left[  \Phi\left(  \mathbf{0}\right)
,\Phi\left(  \mathbf{0}\right)  ,\Phi\left(  \mathbf{1}\right)  \right]   &
=\mathbf{\mu}_{a}^{\left[  3\right]  }\left[  \mathbf{0},\mathbf{0}%
,c\cdot\mathbf{1}\right]  =c\cdot\mathbf{1}.
\end{align}

Thus, the~pair $\left(  \Phi,\Psi\right)  $ is actually the (because~of
(\ref{fa0})) graded homomorphism of the ternary superalgebra $\mathcal{A}%
^{\left[  2,3\right]  }$.

A counterexample to $\Phi$ (\ref{f0}) is $\Phi_{not}\left(  \mathbf{0}\right)
=\mathbf{0+1}$, $\Phi_{not}\left(  \mathbf{1}\right)  =\mathbf{1}$, for~which
$\Phi_{not}\left(  \mathcal{A}\left(  \mathsf{0}\right)  \right)
\nsubseteq\mathcal{A}\left(  \mathsf{0}\right)  $, and~therefore, $\Phi_{not}$
does not preserve~grading.

The same construction can be made for any odd $n$, the~arity of algebra
multiplication, with~only one nonzero product with exactly one
\textquotedblleft fermionic\textquotedblright\ variable $\mathbf{\mu}%
_{a}^{\left[  n\right]  }\left[  \overset{n-1}{\overbrace{\mathbf{0}%
,\ldots,\mathbf{0}}},\mathbf{1}\right]  =\mathbf{1}$ (see (\ref{m0})). This
defines $n$-ary commutative superalgebra $\mathcal{A}^{\left[  2,n\right]  }$,
and $\left(  \Phi,\Psi\right)  $ from (\ref{f0}) and (\ref{f1}) determines its
graded homomorphism.
\end{Example}

We are interested in strictly nonderived groups (or not reducible, meaning
they cannot be reduced to any binary group, even with an automorphism), which
exist for every   $n^{\prime}>2$ \cite{dor3}. In~general, if~$n^{\prime}%
$-ary$\ $group has an idempotent (or identity), it is derived, and~therefore
the first condition to be not strictly derived is%
\begin{equation}
\gcd(N,n^{\prime}-1)>1. \label{gn}%
\end{equation}

This follows from the condition of the neutral element $\mathsf{e}$ in
$\mathsf{G}^{\left[  n^{\prime}\right]  }$ as $\mu_{g}^{\left[  n^{\prime
}\right]  }\left[  \mathsf{e},\ldots\mathsf{e},\mathsf{x}\right]  =\mathsf{x}$
written additively $\left(  n^{\prime}-1\right)  \mathsf{e}+\mathsf{x}%
=\mathsf{x}\operatorname{mod}N$, which gives $\gcd(N,n^{\prime}-1)=1$, or~$N$
and $\left(  n^{\prime}-1\right)  $ are~coprime.

\begin{Example}
[Nonderived ternary grading group]\label{ex-tern}The only strictly nonderived
ternary group with two elements is%
\begin{align}
\mathsf{G}^{\left[  3^{\prime}\right]  }  &  =\left\langle G_{2}\mid\mu
_{g}^{\left[  3^{\prime}\right]  }\right\rangle ,\label{g3}\\
\mu_{g}^{\left[  3^{\prime}\right]  }\left[  \mathsf{x},\mathsf{y}%
,\mathsf{z}\right]   &  =\left(  \mathsf{x}+\mathsf{y}+\mathsf{z}+1\right)
\operatorname{mod}2,\ \ \ \ \mathsf{x},\mathsf{y},\mathsf{z}\in G_{2}.
\label{mx}%
\end{align}

It has $\gcd(2,2)=2$, and~so, no neutral element (identity) exists.
The Cayley table for $\mathsf{G}^{\left[  3\right]  }$ is%
\begin{equation}%
\begin{tabular}
[c]{|c|cc|cc|}\hline
& \ \ \ \ \ $\mathsf{z}$ & $=\mathsf{0}$ & \ \ \ \ \ $\mathsf{z}$ &
$=\mathsf{1}$\\\cline{2-5}
& $\mathsf{y=0}$ & $\mathsf{y=1}$ & $\mathsf{y=0}$ & $\mathsf{y=1}$\\\cline{1-5}
$\mathsf{x=0}$ & $\mathsf{1}$ & \multicolumn{1}{|c|}{$\mathsf{0}$} &
$\mathsf{0}$ & \multicolumn{1}{|c|}{$\mathsf{1}$}\\\cline{1-5}
$\mathsf{x=1}$ & $\mathsf{0}$ & \multicolumn{1}{|c|}{$\mathsf{1}$} &
$\mathsf{1}$ & \multicolumn{1}{|c|}{$\mathsf{0}$}\\\hline
\end{tabular}
\ \label{ct}%
\end{equation}

The corresponding ternary algebra has the decomposition (\ref{a}) in two parts
($\left\vert G\right\vert =N=2$), and~so its addition is binary, because~of
(\ref{sl}) with $m=2$ and $\ell_{m}=1$:%
\begin{equation}
\mathcal{A}^{\left[  2,3\right]  }=\mathcal{A}\left(  \mathsf{0}\right)
\oplus\mathcal{A}\left(  \mathsf{1}\right)  , \label{aa}%
\end{equation}
where $\mathcal{A}\left(  0\right)  $ and $\mathcal{A}\left(  1\right)  $ are
\textquotedblleft even\textquotedblright\ and \textquotedblleft
odd\textquotedblright\ components of the algebra. Using the Cayley
Table~(\ref{ct}), the~conditions (\ref{mas}) in components now become%
\begin{equation}%
\begin{array}
[c]{cc}%
\mathbf{\mu}_{a}^{\left[  3\right]  }\left[  \mathcal{A}\left(  \mathsf{0}%
\right)  ,\mathcal{A}\left(  \mathsf{0}\right)  ,\mathcal{A}\left(
\mathsf{0}\right)  \right]  \subseteq\mathcal{A}\left(  \mathsf{1}\right)
,\ \ \ \ \  & \mathbf{\mu}_{a}^{\left[  3\right]  }\left[  \mathcal{A}\left(
\mathsf{0}\right)  ,\mathcal{A}\left(  \mathsf{0}\right)  ,\mathcal{A}\left(
\mathsf{1}\right)  \right]  \subseteq\mathcal{A}\left(  \mathsf{0}\right)  ,\\
\mathbf{\mu}_{a}^{\left[  3\right]  }\left[  \mathcal{A}\left(  \mathsf{0}%
\right)  ,\mathcal{A}\left(  \mathsf{1}\right)  ,\mathcal{A}\left(
\mathsf{0}\right)  \right]  \subseteq\mathcal{A}\left(  \mathsf{0}\right)
,\ \ \ \ \  & \mathbf{\mu}_{a}^{\left[  3\right]  }\left[  \mathcal{A}\left(
\mathsf{0}\right)  ,\mathcal{A}\left(  \mathsf{1}\right)  ,\mathcal{A}\left(
\mathsf{1}\right)  \right]  \subseteq\mathcal{A}\left(  \mathsf{1}\right)  ,\\
\mathbf{\mu}_{a}^{\left[  3\right]  }\left[  \mathcal{A}\left(  \mathsf{1}%
\right)  ,\mathcal{A}\left(  \mathsf{0}\right)  ,\mathcal{A}\left(
\mathsf{0}\right)  \right]  \subseteq\mathcal{A}\left(  \mathsf{0}\right)
,\ \ \ \ \  & \mathbf{\mu}_{a}^{\left[  3\right]  }\left[  \mathcal{A}\left(
\mathsf{1}\right)  ,\mathcal{A}\left(  \mathsf{0}\right)  ,\mathcal{A}\left(
\mathsf{1}\right)  \right]  \subseteq\mathcal{A}\left(  \mathsf{1}\right)  ,\\
\mathbf{\mu}_{a}^{\left[  3\right]  }\left[  \mathcal{A}\left(  \mathsf{1}%
\right)  ,\mathcal{A}\left(  \mathsf{1}\right)  ,\mathcal{A}\left(
\mathsf{0}\right)  \right]  \subseteq\mathcal{A}\left(  \mathsf{1}\right)
,\ \ \ \ \  & \mathbf{\mu}_{a}^{\left[  3\right]  }\left[  \mathcal{A}\left(
\mathsf{1}\right)  ,\mathcal{A}\left(  \mathsf{1}\right)  ,\mathcal{A}\left(
\mathsf{1}\right)  \right]  \subseteq\mathcal{A}\left(  \mathsf{0}\right)  ,
\end{array}
\label{aaa}%
\end{equation}
{where $\mathbf{\mu}_{a}^{\left[  3\right]  }\mathbf{\ }$is the ternary algebra
multiplication, which is strictly nonderived and can be noncommutative.}

Thus, the~algebra $\mathcal{A}^{\left[  2,3\right]  }$ (\ref{aa}) can be
called a ternary ($3$-ary) superalgebra (even/odd parts in (\ref{aa})) with
ternary ($3^{\prime}$-ary) grading. It can be compared with the ordinary
(binary) superalgebra ($\mathbb{Z}_{2}$-graded algebra), which is in our
notation has the decomposition $\mathcal{A}^{\left[  2,2\right]
}=\left\langle A\mid\mathbf{\nu}_{a}^{\left[  2\right]  },\mathbf{\mu}%
_{a}^{\left[  2\right]  }\right\rangle =\mathcal{A}\left(  \mathsf{0}\right)
\oplus\mathcal{A}\left(  \mathsf{1}\right)  $ with the compatibility
conditions%
\begin{equation}%
\begin{array}
[c]{c}%
\mathbf{\mu}_{a}^{\left[  2\right]  }\left[  \mathcal{A}\left(  \mathsf{0}%
\right)  ,\mathcal{A}\left(  \mathsf{0}\right)  \right]  \subseteq
\mathcal{A}\left(  \mathsf{0}\right)  \mathsf{,}\\
\mathbf{\mu}_{a}^{\left[  2\right]  }\left[  \mathcal{A}\left(  \mathsf{0}%
\right)  ,\mathcal{A}\left(  \mathsf{1}\right)  \right]  \subseteq
\mathcal{A}\left(  1\right)  \mathsf{,\ }\\
\mathbf{\mu}_{a}^{\left[  2\right]  }\left[  \mathcal{A}\left(  \mathsf{1}%
\right)  ,\mathcal{A}\left(  \mathsf{0}\right)  \right]  \subseteq
\mathcal{A}\left(  \mathsf{1}\right)  ,\\
\mathbf{\mu}_{a}^{\left[  2\right]  }\left[  \mathcal{A}\left(  \mathsf{1}%
\right)  ,\mathcal{A}\left(  \mathsf{1}\right)  \right]  \subseteq
\mathcal{A}\left(  \mathsf{0}\right)  .
\end{array}
\label{a0}%
\end{equation}

\end{Example}

\begin{Remark}
In general, the~affine ternary multiplication $\mu_{g}^{\left[  3^{\prime
}\right]  }$ of the $3^{\prime}$-ary grading group $\mathsf{G}^{\left[
3^{\prime}\right]  }$ can be generalized to any $N$ and $n^{\prime}$ as
follows:%
\begin{equation}
\mu_{g}^{\left[  n^{\prime}\right]  }\left[  \mathsf{x}_{1},\mathsf{x}%
_{2},\ldots,\mathsf{x}_{n^{\prime}}\right]  =\left(  \mathsf{x}_{1}%
+\mathsf{x}_{2}+\ldots+\mathsf{x}_{n^{\prime}}+1\right)  \operatorname{mod}%
N,\ \ \ \ \ \mathsf{x}_{1},\mathsf{x}_{2},\ldots,\mathsf{x}_{n^{\prime}}\in
G_{N}, \label{mg}%
\end{equation}
which gives strictly the nonderived grading $n^{\prime}$-ary group, if~the
condition of absence of idempotents (\ref{gn}) holds valid. For~instance, if~
$N=2$, then $n^{\prime}$ should be odd, and the case of $n^{\prime}=N+1$ works always.
\end{Remark}

\begin{Remark}
For higher values of $N$ and $n^{\prime}$, there exist strictly nonderived
grading $n^{\prime}$-ary groups whose multiplication law differs from the
affine form given in (\ref{mg}).
\end{Remark}

\section{\boldmath{$\mathbb{Z}^{\left[  m^{\prime\prime},n^{\prime
\prime}\right]  }$}\textsc{-Graded Polynomials Over }$n${-ary Matrices} \label{sec-poly}}

The next example comes from polynomial rings $\Bbbk\left[  x\right]  $ in one
indeterminate $x$ over a field $\Bbbk$ graded by monomial degree $d$, such
that the direct sum decomposition (\ref{a}) is $\mathcal{R}=\Bbbk\left[
x\right]  =\bigoplus_{d\geq0}\mathcal{R}\left(  d\right)  $, where the nonzero
components are $\mathcal{R}\left(  d\right)  =c\cdot x^{d}$, $c\in\Bbbk$,
$d\in\mathbb{Z}$ (for $d<0$, it is assumed $c=0$). The~elements from
$\mathcal{R}\left(  d\right)  $ are homogeneous of degree $d$, since
$\mathcal{R}\left(  d_{1}\right)  \cap\mathcal{R}\left(  d_{2}\right)
=\varnothing$, and~the polynomial product is%
\begin{equation}
\mathcal{R}\left(  d_{1}\right)  \cdot\mathcal{R}\left(  d_{2}\right)
\subseteq\mathcal{R}\left(  d_{1}+d_{2}\right)  , \label{r}%
\end{equation}
which means that the binary grading group $\mathsf{G}^{\left[  2^{\prime
}\right]  }$ now is the additive group of the ring of integers $\mathbb{Z}$.

The polyadization of this grading group can be done by considering the
additive group of the $\left(  m^{\prime\prime},n^{\prime\prime}\right)
$-ring of polyadic integers $\mathbb{Z}^{\left[  m^{\prime\prime}%
,n^{\prime\prime}\right]  }\left(  a,b\right)  $ consisting of the
representatives of the congruence class $\left[  \left[  a\right]  \right]
_{b}$, ($a,b\in\mathbb{Z}_{+})$, which was introduced in~\cite{dup2017a,duplij2022}. Instead of an indeterminate $x$ we will consider the
block-shift matrices~\cite{dup2022a,duplij2022} which obey $n$-ary
multiplication and binary~addition.

Let us introduce $n$-ary \textquotedblleft matrix
indeterminates\textquotedblright\ as the block-shift matrices of $\left(
n-1\right)  \times\left(  n-1\right)  $ size%
\begin{equation}
X=X\left(  x\right)  =X^{\left[  n\right]  }\left(  x\right)  =M_{\left(
n-1\right)  \times\left(  n-1\right)  }^{\left[  n\right]  }\left(  x\right)
=\left(
\begin{array}
[c]{ccccc}%
0 & x & \ldots & 0 & 0\\
0 & 0 & x & \ldots & 0\\
0 & 0 & \ddots & \ddots & \vdots\\
\vdots & \vdots & \ddots & 0 & x\\
x & 0 & \ldots & 0 & 0
\end{array}
\right)  , \label{mb}%
\end{equation}
which are closed under $n$-ary multiplication, which is strictly nonderived~\cite{dup2022a}. Without~any additional requirements the set of $n$-ary
matrices (\ref{mb}) form the totally commutative   $n$-ary semigroup%
\begin{align}
\mathcal{M}^{\left[  n\right]  }  &  =\left\langle \left\{  X\right\}  \mid
\mu_{x}^{\left[  n\right]  }\right\rangle ,\label{mnx}\\
\mu_{x}^{\left[  n\right]  }\left[  X\left(  x_{1}\right)  ,X\left(
x_{2}\right)  ,\ldots,X\left(  x_{n}\right)  \right]   &  =X\left(
x_{1}\right)  \cdot X\left(  x_{2}\right)  \cdot\ldots\cdot X\left(
x_{n}\right)  , \label{mn1}%
\end{align}
where the matrix product in the r.h.s. is in $\Bbbk$. Obviously,
\begin{equation}
\mu_{x}^{\left[  n\right]  }\left[  X\left(  x_{1}\right)  ,X\left(
x_{2}\right)  ,\ldots,X\left(  x_{n}\right)  \right]  =X\left(  x_{1}%
x_{2}\ldots x_{n}\right)  .
\end{equation}

If $x$ is invertible, then $\mathcal{M}^{\left[  n\right]  }$ becomes an
$n$-ary group with the querelement%
\begin{align}
\bar{X}  &  =\bar{X}\left(  x\right)  =X\left(  x^{2-n}\right)  \left(
\begin{array}
[c]{ccccc}%
0 & x^{2-n} & \ldots & 0 & 0\\
0 & 0 & x^{2-n} & \ldots & 0\\
0 & 0 & \ddots & \ddots & \vdots\\
\vdots & \vdots & \ddots & 0 & x^{2-n}\\
x^{2-n} & 0 & \ldots & 0 & 0
\end{array}
\right)  ,\label{xx}\\
\mu_{x}^{\left[  n\right]  }\left[  \overset{n-1}{\overbrace{X,X,\ldots,X}%
},\bar{X}\right]   &  =X,\ \ n\geq3. \label{mxx}%
\end{align}

The polyadic identity in $\mathcal{M}^{\left[  n\right]  }$ is%
\begin{equation}
E^{\left[  n\right]  }=X\left(  x^{0}\right)  =X\left(  1\right)  ,
\end{equation}
as%
%\vspace{-15pt}
%\begin{adjustwidth}{-\extralength}{0cm}
\begin{align}
E^{\left[  n\right]  }  &  =\left(
\begin{array}
[c]{ccccc}%
0 & 1 & \ldots & 0 & 0\\
0 & 0 & 1 & \ldots & 0\\
0 & 0 & \ddots & \ddots & \vdots\\
\vdots & \vdots & \ddots & 0 & 1\\
1 & 0 & \ldots & 0 & 0
\end{array}
\right)  _{\left(  n-1\right)  \times\left(  n-1\right)  }\neq I^{\left[
n\right]  }=\left(
\begin{array}
[c]{ccccc}%
1 & 0 & \ldots & 0 & 0\\
0 & 1 & 0 & \ldots & 0\\
0 & 0 & \ddots & \ddots & \vdots\\
\vdots & \vdots & \ddots & 1 & 0\\
0 & 0 & \ldots & 0 & 1
\end{array}
\right)  _{\left(  n-1\right)  \times\left(  n-1\right)  }\notin
\mathcal{M}^{\left[  n\right]  },\ n\geq3,\label{e}\\
E^{\left[  2\right]  }  &  =I^{\left[  2\right]  }=\left(
\begin{array}
[c]{cc}%
1 & 0\\
0 & 1
\end{array}
\right)  ,\label{ei}\\
&  \mu_{x}^{\left[  n\right]  }\left[  \overset{n-1}{\overbrace{E^{\left[
n\right]  },E^{\left[  n\right]  },\ldots,E^{\left[  n\right]  }}},X\left(
x\right)  \right]  =X\left(  x\right)  ,\ \ n\geq2. \label{me}%
\end{align}
%\end{adjustwidth}

Due to the \textquotedblleft quantization\textquotedblright\ (\ref{w}), we can
multiply only the admissible number of the matrix indeterminates, and the
degree (word length of $x$'s) is also \textquotedblleft
quantized\textquotedblright:%
\begin{equation}
d_{\ell}=d_{\ell}^{\left[  n\right]  }=\ell\left(  n-1\right)  +1, \label{dl}%
\end{equation}
where $\ell\in\mathbb{N}_{0}$ is the polyadic power (\ref{xl}). So, we have the
polyadic monomials as polyadic powers (instead of $x^{d}$)%
\begin{align}
X^{\left\langle \ell\right\rangle }  &  \equiv X^{d_{\ell}}=X\left(
x^{d_{\ell}}\right)  =M_{\left(  n-1\right)  \times\left(  n-1\right)
}^{\left[  n\right]  }\left(  x^{d_{\ell}}\right)  ,\label{xm}\\
X^{0}  &  =I^{\left[  n\right]  }\notin\mathcal{M}^{\left[  n\right]  }.
\label{x0}%
\end{align}

\begin{Remark}
\label{rem-free}In the binary case $n=2$, the degree (\ref{dl}) becomes
$d_{\ell}=\ell+1$, and~so, the zero polyadic power corresponds to one degree. Let
informally $\ell=-1$; then, in the general case, (\ref{dl}) is $d_{-1}=2-n$.
Thus, the free term (having zero degree $d_{\ell}=0$) of a polynomial belongs
to $\mathcal{M}^{\left[  n\right]  }$ only in the binary case; see (\ref{e})
and (\ref{x0}). Note that negative polyadic power is calculated for $n$-ary
groups in the special way; that is, $X^{\left\langle -1\right\rangle }%
\equiv\bar{X}$ \cite{duplij2022}. Therefore, the~free term in the $n$-ary case
becomes $c_{0}\cdot X\left(  x^{0}\right)  =c_{0}\cdot E^{\left[  n\right]  }$
(as for ordinary polynomial $c_{0}\cdot x^{0}=c_{0}\cdot1=c_{0}$).
\end{Remark}

Thus, the~generic polynomial over $n$-ary matrices of the length $L+1$ with
indeterminate $x$ has the form%
%\vspace{-12pt}
%\begin{adjustwidth}{-\extralength}{0cm}
\begin{align}
P^{\left[  n\right]  }  &  =P^{\left[  n\right]  }\left(  x\right)  =c\cdot
E^{\left[  n\right]  }+\sum_{\ell=0}^{L}c_{\ell}\cdot X^{\left[  n\right]
}\left(  x^{d_{\ell}}\right)  =c\cdot E^{\left[  n\right]  }+c_{0}\cdot
X^{\left[  n\right]  }\left(  x\right)  +c_{1}\cdot X^{\left[  n\right]
}\left(  x^{n}\right) \nonumber\\
&  +c_{2}\cdot X^{\left[  n\right]  }\left(  x^{2n-1}\right)  +\ldots c_{\ell
}\cdot X^{\left[  n\right]  }\left(  x^{\ell\left(  n-1\right)  +1}\right)
+\ldots+c_{L}\cdot X^{\left[  n\right]  }\left(  x^{L\left(  n-1\right)
+1}\right)  ,\ \ c_{\ell}\in\Bbbk. \label{pn}%
\end{align}
%\end{adjustwidth}

As for ordinary polynomials, $c_{i}\cdot X^{\left[  n\right]  }\left(
x^{i}\right)  \cap c_{j}\cdot X^{\left[  n\right]  }\left(  x^{j}\right)
=\varnothing$, if~$i\neq j$, and~therefore, the sum (\ref{pn}) is~direct.

{The polyadic matrix algebra corresponding to $\Bbbk\left[  X\right]  $ has the
similar direct sum decomposition}
\begin{equation}
\mathcal{A}^{\left[  2,n\right]  }=\mathcal{A}_{0}\oplus\bigoplus_{\ell
=0}\mathcal{A}\left(  d_{\ell}\right)  , \label{adl}%
\end{equation}
where $\mathcal{A}_{0}=c\cdot E^{\left[  n\right]  }$, $\mathcal{A}\left(
d_{\ell}\right)  \mathcal{=}c\cdot M_{\left(  n-1\right)  \times\left(
n-1\right)  }^{\left[  n\right]  }\left(  x^{d_{\ell}}\right)  =c_{\ell}\cdot
X\left(  x^{d_{\ell}}\right)  $l see (\ref{xm}) and (\ref{pn}). Obviously,
$\mathcal{A}\left(  d_{\ell_{1}}\right)  \cap\mathcal{A}\left(
d_{\ell_{2}}\right)  =\varnothing$, if~$l_{1}\neq l_{2}$, as~in the binary
case, and~the $n$-ary multiplication is given by the ordinary matrix product in
$\Bbbk$:%

\begin{equation}
\mathbf{\mu}_{a}^{\left[  n\right]  }\left[  \mathbf{a}_{1},\mathbf{a}%
_{2},\ldots,\mathbf{a}_{n}\right]  =\mathbf{a}_{1}\cdot\mathbf{a}_{2}%
\cdot\ldots\cdot\mathbf{a}_{n},\ \ \ \ \mathbf{a}_{j}\in\mathcal{A}^{\left[
2,n\right]  }.
\end{equation}

Then, the $n$-ary product of components (\ref{ma}) becomes%
\begin{equation}
\mathbf{\mu}_{a}^{\left[  n\right]  }\left[  \mathcal{A}\left(  d_{\ell_{1}%
}\right)  ,\mathcal{A}\left(  d_{\ell_{2}}\right)  ,\ldots,\mathcal{A}\left(
d_{\ell_{n}}\right)  \right]  \subseteq\mathcal{A}\left(  d_{\ell}\right)
,\ \ \ \ \ \ell=\ell_{1}+\ell_{2}+\ldots+\ell_{n}. \label{mad}%
\end{equation}

Recall that the polyadic integers~\cite{dup2017a} are special representatives
$y_{k}=a+kb$ of a congruence class $\left[  \left[  a\right]  \right]
_{b}=\left\{  y_{k}\right\}  $, $k\in\mathbb{Z}$, $b>0$, $0\leq a\leq b-1$,
$a,b\in\mathbb{Z}_{+}$, which satisfy the \textquotedblleft
quantization\textquotedblright\ conditions (the arity shape invariants should
be positive integers)%
\begin{align}
I^{\left[  m^{\prime\prime}\right]  }\left(  a,b\right)   &  =\frac{a}%
{b}\left(  m^{\prime\prime}-1\right)  \in\mathbb{Z}_{+},\label{i}\\
J^{\left[  n^{\prime\prime}\right]  }\left(  a,b\right)   &  =\frac{a}%
{b}\left(  a^{n^{\prime\prime}-1}-1\right)  \in\mathbb{Z}_{+}. \label{j}%
\end{align}
In this case only, the~representatives form the strictly nonderived $\left(
m^{\prime\prime},n^{\prime\prime}\right)  $-ring%
\begin{align}
\mathbb{Z}^{\left[  m^{\prime\prime},n^{\prime\prime}\right]  }\left(
a,b\right)   &  =\left\langle \left\{  y_{k}=a+kb\right\}  \mid\nu^{\left[
m^{\prime\prime}\right]  },\mu^{\left[  n^{\prime\prime}\right]
}\right\rangle ,\label{z}\\
\nu^{\left[  m^{\prime\prime}\right]  }\left[  y_{1},y_{2},\ldots
,y_{m^{\prime\prime}}\right]   &  =y_{1}+y_{2}+\ldots+y_{m^{\prime\prime}},\\
\mu^{\left[  n^{\prime\prime}\right]  }\left[  y_{1},y_{2},\ldots
,y_{n^{\prime\prime}}\right]   &  =y_{1}\cdot y_{2}\cdot\ldots\cdot
y_{n^{\prime\prime}},
\end{align}
which is a polyadic analog of the ring of ordinary integers $\mathbb{Z=Z}%
^{\left[  2,2\right]  }\left(  0,1\right)  $, for~which both shape invariants
(\ref{i}) {and} %MDPI: We have revised the en dash to ``and'' please check and confirm. The~following highlights are the same
%%%%%%Author's ANSWER: Confirmed everywhere.
 (\ref{j}) vanish:%
\begin{equation}
I^{\left[  m^{\prime\prime}\right]  }\left(  a,b\right)  =0,\ J^{\left[
n^{\prime\prime}\right]  }\left(  a,b\right)  =0\Longleftrightarrow
\mathbb{Z}^{\left[  m^{\prime\prime},n^{\prime\prime}\right]  }\left(
a,b\right)  =\mathbb{Z}.
\end{equation}

Now, we use the additive $m^{\prime\prime}$-ary group of $\mathbb{Z}^{\left[
m^{\prime\prime},n^{\prime\prime}\right]  }\left(  a,b\right)  $ as the
grading group $\mathsf{G}^{\left[  n^{\prime}\right]  }$ for polyadic
polynomials, by~full analogy with the binary case (\ref{r}). This means that
the consistency condition of the grading (\ref{ma}) should be%
\begin{equation}
\mathbf{\mu}_{a}^{\left[  n\right]  }\left[  \mathcal{A}\left(  d_{\ell_{1}%
}\right)  ,\mathcal{A}\left(  d_{\ell_{2}}\right)  ,\ldots,\mathcal{A}\left(
d_{\ell_{n}}\right)  \right]  \subseteq\mathcal{A}\left(  \nu^{\left[
m^{\prime\prime}\right]  }\left[  y_{1},y_{2},\ldots,y_{m^{\prime\prime}%
}\right]  \right)  . \label{ay}%
\end{equation}

It follows from (\ref{ay}) that the arities of algebra $\mathcal{A}^{\left[
2,n\right]  }$ multiplication and the grading ring $\mathbb{Z}^{\left[
m^{\prime\prime},n^{\prime\prime}\right]  }\left(  a,b\right)  $ addition
(being the arity of the grading group $\mathsf{G}^{\left[  n^{\prime}\right]
}$) coincide as%
\begin{equation}
n^{\prime}=m^{\prime\prime}=n. \label{nmn}%
\end{equation}

We assume that every degree of the monomial $c_{k}\cdot X^{d_{l_{k}}}$ in
$\mathcal{A}^{\left[  2,n\right]  }$ is equal to the corresponding
representative of the grading polyadic ring in (\ref{ay}), which gives the
relations (together with (\ref{nmn}))%
\begin{align}
d_{\ell_{k}}  &  =y_{k},\label{dy}\\
\ell_{k}\left(  n-1\right)  +1  &  =a+b\cdot k. \label{la}%
\end{align}

One set of solutions to (\ref{la}) is given by ($n\geq3$)%
\begin{align}
a  &  =1,\label{q1}\\
b  &  =n-1,\label{b1}\\
k  &  =\ell_{k}. \label{k}%
\end{align}

This means that a polyadic algebra with given arity of multiplication $n$ can be
graded not by arbitrary polyadic integers, as~in the binary case, but by
those coming from the special polyadic rings $\mathbb{Z}^{\left[
n,n^{\prime\prime}\right]  }\left(  1,n-1\right)  $ having the same arity of
addition $m^{\prime\prime}=n$, arbitrary arity of multiplication
$n^{\prime\prime}$ and the fixed arity shape invariants (\ref{i}) and (\ref{j}):%
\begin{align}
I^{\left[  n\right]  }\left(  1,n-1\right)   &  =1,\label{i1}\\
J^{\left[  n^{\prime\prime}\right]  }\left(  1,n-1\right)   &  =0. \label{j1}%
\end{align}

\begin{Example}
Let us consider the concrete polynomial over $4$-ary (block-shift) matrices
with the indeterminate $x$, as%
\begin{align}
X  &  =X^{\left[  4\right]  }\left(  x\right)  =M_{3\times3}^{\left[
4\right]  }\left(  x\right)  =\left(
\begin{array}
[c]{ccc}%
0 & x & 0\\
0 & 0 & x\\
x & 0 & 0
\end{array}
\right)  ,\ \ \ E=E^{\left[  6\right]  }=\left(
\begin{array}
[c]{ccc}%
0 & 1 & 0\\
0 & 0 & 1\\
1 & 0 & 0
\end{array}
\right)  ,\label{x6}\\
P  &  =P^{\left[  4\right]  }=3E-12X^{7}+7X^{10}+5X^{16}-8X^{19}. \label{p6}%
\end{align}

The degree (\ref{dl}) now becomes%
\begin{equation}
d_{\ell}=d_{\ell}^{\left[  4\right]  }=3\ell+1. \label{d3}%
\end{equation}

The monomials are as follows: free constant term $3E\in\mathcal{A}_{0}$, $\left(
-12X^{7}\right)  =\left(  -12X^{d_{2}}\right)  \in\mathcal{A}\left(
10\right)  =\mathcal{A}\left(  d_{2}\right)  $ has polyadic power $\ell=2$,
$7X^{10}=7X^{d_{3}}\in\mathcal{A}\left(  10\right)  =\mathcal{A}\left(
d_{3}\right)  $ has polyadic power $\ell=3$, $5X^{16}=5X^{d_{5}}\in
\mathcal{A}\left(  19\right)  =\mathcal{A}\left(  d_{5}\right)  $ has the
polyadic power $\ell=5$, and $\left(  -8X^{19}\right)  =\left(  -8X^{d_{6}%
}\right)  \in\mathcal{A}\left(  19\right)  $ has the polyadic power $\ell=6$.
The sum (\ref{p6}) is direct, evidently, and~so, we obtain the unique direct
sum decomposition of the corresponding (to $\Bbbk\left[  X\right]  $) algebra%
%\vspace{-8pt}
%\begin{adjustwidth}{-\extralength}{0cm}
\begin{equation}
\mathcal{A}^{\left[  2,4\right]  }=\mathcal{A}_{0}\oplus\mathcal{A}\left(
7\right)  \oplus\mathcal{A}\left(  10\right)  \oplus\mathcal{A}\left(
16\right)  \oplus\mathcal{A}\left(  19\right)  =\mathcal{A}_{0}\oplus
\mathcal{A}\left(  d_{2}\right)  \oplus\mathcal{A}\left(  d_{3}\right)
\oplus\mathcal{A}\left(  d_{5}\right)  \oplus\mathcal{A}\left(  d_{6}\right)
. \label{a24}%
\end{equation}
%\end{adjustwidth}

Then, we choose the grading polyadic ring (\ref{z}) using (\ref{q1})--(\ref{k})
for $a=1$ and $b=3$ and~take, for~instance, $\mathbb{Z}^{\left[  4,7\right]
}\left(  1,3\right)  $. The~additive $4$-ary group of this ring should be used
to \textquotedblleft add\textquotedblright\ the degrees, while multiplying $4$
polynomials of the fixed degrees in (\ref{ay}). For~instance, the~main
consistency condition (\ref{ay}) for all nonconstant graded components of the
polynomial $p^{\left[  4\right]  }$ (\ref{p6}) has the form (using (\ref{la}))%
%\vspace{-15pt}
%\begin{adjustwidth}{-\extralength}{0cm}
\begin{align}
\mathbf{\mu}_{a}^{\left[  4\right]  }\left[  \mathcal{A}\left(  d_{2}\right)
,\mathcal{A}\left(  d_{3}\right)  ,\mathcal{A}\left(  d_{5}\right)
,\mathcal{A}\left(  d_{6}\right)  \right]   &  \subseteq\mathcal{A}\left(
\nu^{\left[  4\right]  }\left[  \left(  1+3\cdot2\right)  ,\left(
1+3\cdot3\right)  ,\left(  1+3\cdot5\right)  ,\left(  1+3\cdot6\right)
\right]  \right) \nonumber\\
&  =\mathcal{A}\left(  1+3\cdot17\right)  =\mathcal{A}\left(  52\right)  ,
\end{align}
%\end{adjustwidth}
where on r.h.s. there are elements of the polyadic ring $\mathbb{Z}^{\left[
4,7\right]  }\left(  1,3\right)  $.

Thus, the~multiplication in the polyadic algebra $\mathcal{A}^{\left[
2,4\right]  }$ (\ref{a24}), corresponding to the polynomial (\ref{p6}) over
$4$-ary matrices, respects the multiary grading by the polyadic ring
$\mathbb{Z}^{\left[  4,7\right]  }\left(  1,3\right)  $.
\end{Example}

\section{\label{sec-high}\textsc{Higher Power Multiary~Gradings}}

The compatibility condition (\ref{ma}) in the higher arity case can have a more
complicated structure, leading to inequality $n^{\prime}\neq n$ instead of
(\ref{mn}), if~we use polyadic powers and the \textquotedblleft
quantization\textquotedblright\ of word length (\ref{w}).

\begin{Definition}
A multiary higher power $\mathsf{G}$-graded polyadic $\Bbbk$-algebra is the
direct sum decomposition of $\Bbbk$-vector spaces (\ref{a}) such that the
$n$-ary multiplication in the algebra respects the $n^{\prime}$-ary
multiplication in the grading group as follows:%
\begin{equation}
\mathbf{\mu}_{a}^{\left[  n\right]  \circ\ell_{n}}\left[  \mathcal{A}\left(
\mathsf{g}_{1}\right)  ,\mathcal{A}\left(  \mathsf{g}_{2}\right)
,\ldots,\mathcal{A}\left(  \mathsf{g}_{\ell_{n}\left(  n-1\right)  +1}\right)
\right]  \subseteq\mathcal{A}\left(  \mathbf{\mu}_{g}^{\left[  n^{\prime
}\right]  \circ\ell_{n^{\prime}}}\left[  \mathsf{g}_{1},\mathsf{g}_{2}%
,\ldots,\mathsf{g}_{\ell_{n^{\prime}}\left(  n^{\prime}-1\right)  +1}\right]
\right)  , \label{ma1}%
\end{equation}
where $\mathcal{A}\left(  \mathsf{g}_{i}\right)  $ is the $i$th component of
the decomposition (\ref{a}), while $\ell_{n}$ and $\ell_{n^{\prime}}$ are
polyadic powers (\ref{xl}) of the algebra and the grading group
multiplications, respectively.

A higher power polyadic algebra is strongly graded, if,~in (\ref{ma1}), the
equality is%
\begin{equation}
\mathbf{\mu}_{a}^{\left[  n\right]  \circ\ell_{n}}\left[  \mathcal{A}\left(
\mathsf{g}_{1}\right)  ,\mathcal{A}\left(  \mathsf{g}_{2}\right)
,\ldots,\mathcal{A}\left(  \mathsf{g}_{\ell_{n}\left(  n-1\right)  +1}\right)
\right]  =\mathcal{A}\left(  \mathbf{\mu}_{g}^{\left[  n^{\prime}\right]
\circ\ell_{n^{\prime}}}\left[  \mathsf{g}_{1},\mathsf{g}_{2},\ldots
,\mathsf{g}_{\ell_{n^{\prime}}\left(  n^{\prime}-1\right)  +1}\right]
\right)  . \label{mas1}%
\end{equation}

\end{Definition}

\begin{Theorem}
\label{theor-l=l}The arity of multiplication of the higher power polyadic
algebra $\mathcal{A}^{\left[  m,n\right]  }$ and the arity of the multiary
grading group $\mathsf{G}^{\left[  n^{\prime}\right]  }$ are \textquotedblleft
quantized\textquotedblright\ and connected by (cf. (\ref{mn})), as follows:
\begin{equation}
\ell_{n^{\prime}}\left(  n^{\prime}-1\right)  =\ell_{n}\left(  n-1\right)  .
\label{ll}%
\end{equation}

\end{Theorem}

\begin{proof}
The definition of the polyadic power
as the number of operations in their composition (\ref{w}) follows from (\ref{ma1}) {and} (\ref{mas1}), and~then, by equating
the allowed word lengths,   $w=\ell_{n^{\prime}}\left(  n^{\prime}-1\right)
+1=\ell_{n}\left(  n-1\right)  +1$ in both sides.
\end{proof}

\begin{center}
The solutions of (\ref{ll}) for $n^{\prime}\neq n$, which are $\leq$ $5$
together with the word length $w$, are %
\begin{equation}%
\begin{tabular}
[c]{|c|c|c|c|c|}\hline
$\ell_{n^{\prime}}$ & $\ell_{n}$ & $n^{\prime}$ & $n$ & $w$\\\hline
3 & 2 & 3 & 4 & \textbf{7}\\
2 & 1 & 3 & 5 & \textbf{5}\\
4 & 2 & 3 & 5 & \textbf{9}\\
2 & 3 & 4 & 3 & \textbf{7}\\
4 & 3 & 4 & 5 & \textbf{13}\\
1 & 2 & 5 & 3 & \textbf{5}\\
2 & 4 & 5 & 3 & \textbf{9}\\
3 & 4 & 5 & 4 & \textbf{13}\\\hline
\end{tabular}
\ \label{ln}%
\end{equation}

\end{center}

The \textquotedblleft quantization\textquotedblright\ condition (\ref{ll})
together with {Table} %MDPI: We revised ``Table'' to ``Equation''. Please confirm.
%%%%%%Author's ANSWER: Better to retain `Table', because (124) _is_ the table!
~(\ref{ln}) shows that multiary gradings of higher arity
polyadic algebras are possible only for certain specific combinations of
arities and powers. In~contrast, the~binary case imposes no such~restrictions.

\begin{Example}
Let us consider the ternary grading group from \textit{Example} \ref{ex-tern}
as in (\ref{g3}) and (\ref{mx}). To~satisfy (\ref{ll}), we take the second solution
in (\ref{ln}) $\ell_{n^{\prime}}=2$, $\ell_{n}=1$, $n^{\prime}=3$, $n=5$,
which shows that the higher power polyadic algebra with binary addition
(because of $m=2$ in the decomposition (\ref{aa})) should have $5$-ary
multiplication, and~so, it is $\mathcal{A}^{\left[  2,5\right]  }$. Since
$\ell_{n^{\prime}}=2$, the~higher (third) polyadic power $3^{\prime}$-ary
grading group multiplication becomes%
\begin{equation}
\mathbf{\mu}_{g}^{\left[  3^{\prime}\right]  \circ2}\left[  \mathsf{x}%
,\mathsf{y},\mathsf{z},\mathsf{t},\mathsf{u}\right]  =\left(  \mathsf{x}%
+\mathsf{y}+\mathsf{z}+\mathsf{t}+\mathsf{u}+1\right)  \operatorname{mod}%
2,\ \ \ \ \mathsf{x},\mathsf{y},\mathsf{z},\mathsf{t},\mathsf{u}\in G_{2}.
\label{m33}%
\end{equation}

The Cayley table of (\ref{m33}) is%
\begin{equation}
{\tiny
\begin{tabular}
[c]{|crc|rc|}\hline
$(\mathsf{z},\mathsf{t},\mathsf{u})\in$ & $\left\{
\mathsf{(0,0,0),(1,1,0),(1,0,1),(0,1,1)}\right\}  $ &  & $\left\{
\mathsf{(1,0,0),(0,1,0),(0,0,1),(1,1,1)}\right\}  $ & \\\hline
& \multicolumn{1}{|c}{$\mathsf{y=0}$} & $\mathsf{y=1}$ &
\multicolumn{1}{|c}{$\mathsf{y=0}$} & $\mathsf{y=1}$\\\hline
$\mathsf{x=0}$ & \multicolumn{1}{|c}{$\mathsf{1}$} &
\multicolumn{1}{|c|}{$\mathsf{0}$} & \multicolumn{1}{|c}{$\mathsf{0}$} &
\multicolumn{1}{|c|}{$\mathsf{1}$}\\\hline
$\mathsf{x=1}$ & \multicolumn{1}{|c}{$\mathsf{0}$} &
\multicolumn{1}{|c|}{$\mathsf{1}$} & \multicolumn{1}{|c}{$\mathsf{1}$} &
\multicolumn{1}{|c|}{$\mathsf{0}$}\\\hline
\end{tabular}
\ \ \label{ct5}}%
\end{equation}

The corresponding polyadic algebra has the decomposition (\ref{a}) in two
parts ($\left\vert G\right\vert =N=2$):%
\begin{equation}
\mathcal{A}^{\left[  2,5\right]  }=\mathcal{A}\left(  \mathsf{0}\right)
\oplus\mathcal{A}\left(  \mathsf{1}\right)  , \label{a25}%
\end{equation}
where $\mathcal{A}\left(  \mathsf{0}\right)  $ and $\mathcal{A}\left(
\mathsf{1}\right)  $ are the \textquotedblleft even\textquotedblright\ and
\textquotedblleft odd\textquotedblright\ components of the algebra. Using the
Cayley  {table}~(\ref{ct5}), the~conditions (\ref{mas}) in components now become%
%%%%%%Author's ANSWER: Better to retain `Table', because (126) _is_ the table! There is no ``Cayley equation'' in math science!
%\vspace{-5pt}
%\begin{adjustwidth}{-\extralength}{0cm}
\[%
\begin{array}
[c]{ll}%
\mathbf{\mu}_{a}^{\left[  5\right]  }\left[  \mathcal{A}\left(  0\right)
,\mathcal{A}\left(  0\right)  ,\mathcal{A}\left(  0\right)  ,\mathcal{A}%
\left(  0\right)  ,\mathcal{A}\left(  0\right)  \right]  \subseteq
\mathcal{A}\left(  1\right)  , & \mathbf{\mu}_{a}^{\left[  5\right]  }\left[
\mathcal{A}\left(  0\right)  ,\mathcal{A}\left(  0\right)  ,\mathcal{A}\left(
1\right)  ,\mathcal{A}\left(  0\right)  ,\mathcal{A}\left(  1\right)  \right]
\subseteq\mathcal{A}\left(  1\right)  ,\\
\mathbf{\mu}_{a}^{\left[  5\right]  }\left[  \mathcal{A}\left(  0\right)
,\mathcal{A}\left(  1\right)  ,\mathcal{A}\left(  0\right)  ,\mathcal{A}%
\left(  0\right)  ,\mathcal{A}\left(  0\right)  \right]  \subseteq
\mathcal{A}\left(  0\right)  , & \mathbf{\mu}_{a}^{\left[  5\right]  }\left[
\mathcal{A}\left(  0\right)  ,\mathcal{A}\left(  1\right)  ,\mathcal{A}\left(
1\right)  ,\mathcal{A}\left(  0\right)  ,\mathcal{A}\left(  1\right)  \right]
\subseteq\mathcal{A}\left(  0\right)  ,\\
\mathbf{\mu}_{a}^{\left[  5\right]  }\left[  \mathcal{A}\left(  1\right)
,\mathcal{A}\left(  0\right)  ,\mathcal{A}\left(  0\right)  ,\mathcal{A}%
\left(  0\right)  ,\mathcal{A}\left(  0\right)  \right]  \subseteq
\mathcal{A}\left(  0\right)  , & \mathbf{\mu}_{a}^{\left[  5\right]  }\left[
\mathcal{A}\left(  1\right)  ,\mathcal{A}\left(  0\right)  ,\mathcal{A}\left(
1\right)  ,\mathcal{A}\left(  0\right)  ,\mathcal{A}\left(  1\right)  \right]
\subseteq\mathcal{A}\left(  0\right)  ,\\
\mathbf{\mu}_{a}^{\left[  5\right]  }\left[  \mathcal{A}\left(  1\right)
,\mathcal{A}\left(  1\right)  ,\mathcal{A}\left(  0\right)  ,\mathcal{A}%
\left(  0\right)  ,\mathcal{A}\left(  0\right)  \right]  \subseteq
\mathcal{A}\left(  1\right)  , & \mathbf{\mu}_{a}^{\left[  5\right]  }\left[
\mathcal{A}\left(  1\right)  ,\mathcal{A}\left(  1\right)  ,\mathcal{A}\left(
1\right)  ,\mathcal{A}\left(  0\right)  ,\mathcal{A}\left(  1\right)  \right]
\subseteq\mathcal{A}\left(  1\right)  ,
\end{array}
\]
%
%\end{adjustwidth}

%\begin{adjustwidth}{-\extralength}{0cm}
\begin{equation}%
\begin{array}
[c]{ll}%
\mathbf{\mu}_{a}^{\left[  5\right]  }\left[  \mathcal{A}\left(  0\right)
,\mathcal{A}\left(  0\right)  ,\mathcal{A}\left(  0\right)  ,\mathcal{A}%
\left(  1\right)  ,\mathcal{A}\left(  0\right)  \right]  \subseteq
\mathcal{A}\left(  0\right)  , & \mathbf{\mu}_{a}^{\left[  5\right]  }\left[
\mathcal{A}\left(  0\right)  ,\mathcal{A}\left(  0\right)  ,\mathcal{A}\left(
1\right)  ,\mathcal{A}\left(  1\right)  ,\mathcal{A}\left(  1\right)  \right]
\subseteq\mathcal{A}\left(  0\right)  ,\\
\mathbf{\mu}_{a}^{\left[  5\right]  }\left[  \mathcal{A}\left(  0\right)
,\mathcal{A}\left(  1\right)  ,\mathcal{A}\left(  0\right)  ,\mathcal{A}%
\left(  1\right)  ,\mathcal{A}\left(  0\right)  \right]  \subseteq
\mathcal{A}\left(  1\right)  , & \mathbf{\mu}_{a}^{\left[  5\right]  }\left[
\mathcal{A}\left(  0\right)  ,\mathcal{A}\left(  1\right)  ,\mathcal{A}\left(
1\right)  ,\mathcal{A}\left(  1\right)  ,\mathcal{A}\left(  1\right)  \right]
\subseteq\mathcal{A}\left(  1\right)  ,\\
\mathbf{\mu}_{a}^{\left[  5\right]  }\left[  \mathcal{A}\left(  1\right)
,\mathcal{A}\left(  0\right)  ,\mathcal{A}\left(  0\right)  ,\mathcal{A}%
\left(  1\right)  ,\mathcal{A}\left(  0\right)  \right]  \subseteq
\mathcal{A}\left(  1\right)  , & \mathbf{\mu}_{a}^{\left[  5\right]  }\left[
\mathcal{A}\left(  1\right)  ,\mathcal{A}\left(  0\right)  ,\mathcal{A}\left(
1\right)  ,\mathcal{A}\left(  1\right)  ,\mathcal{A}\left(  1\right)  \right]
\subseteq\mathcal{A}\left(  1\right)  ,\\
\mathbf{\mu}_{a}^{\left[  5\right]  }\left[  \mathcal{A}\left(  1\right)
,\mathcal{A}\left(  1\right)  ,\mathcal{A}\left(  0\right)  ,\mathcal{A}%
\left(  1\right)  ,\mathcal{A}\left(  0\right)  \right]  \subseteq
\mathcal{A}\left(  0\right)  , & \mathbf{\mu}_{a}^{\left[  5\right]  }\left[
\mathcal{A}\left(  1\right)  ,\mathcal{A}\left(  1\right)  ,\mathcal{A}\left(
1\right)  ,\mathcal{A}\left(  1\right)  ,\mathcal{A}\left(  1\right)  \right]
\subseteq\mathcal{A}\left(  0\right)  .
\end{array}
\label{a5}%
\end{equation}
%\end{adjustwidth}

We call the algebra $\mathcal{A}^{\left[  2,5\right]  }$ (\ref{a5}) a $5$-ary
superalgebra (even/odd parts in (\ref{a25})) with higher (second) power
ternary ($3^{\prime}$-ary) grading (\ref{m33}).

Note that this $5$-ary superalgebra (\ref{a5}) cannot be reduced to the
ordinary (binary) superalgebra (\ref{a0}), because~the ternary grading group
$\mathsf{G}^{\left[  3^{\prime}\right]  }$ is not strictly derived (not
reduced to a binary group).
\end{Example}

\section{\textsc{Conclusions}}

In this work, we have established a comprehensive framework for multiary
graded polyadic algebras, extending classical grading theory to higher-arity
algebraic structures. Our investigation has revealed several fundamentally new
phenomena that distinguish the polyadic case from its binary~counterpart.

\textbf{{Main} 
 contributions}

1. \textit{{General} 
 framework}: We introduced the concept of multiary
$\mathsf{G}$-graded polyadic algebras, defined by the decomposition
$\mathcal{A}^{[m,n]}=\bigoplus_{\mathsf{g}_{i}\in\mathsf{G}^{[n^{\prime}]}%
}\mathcal{A}(\mathsf{g}_{i})$ with the compatibility condition $\boldsymbol{\mu
}_{a}^{[n]}[\mathcal{A}(\mathsf{g}_{1}),\ldots,\mathcal{A}(\mathsf{g}%
_{n})]\subseteq\mathcal{A}(\mu_{g}^{[n^{\prime}]}[\mathsf{g}_{1}%
,\ldots,\mathsf{g}_{n^{\prime}}])$. This definition accommodates arbitrary
arities for both algebra operations and grading group~operations.

2. \textit{{Quantization} rules}: We discovered precise constraints connecting
the arities:

\begin{itemize}
\item For strongly $\mathsf{G}$-graded algebras, the arity of graded group
coincides with the multiplication arity of algebra $n^{\prime}=n$
({Proposition} \ref{prop-n=n}).

\item The relation between the grading group order and the arity of the algebra
addition is $|G|=\ell_{m}(m-1)+1$ ({Theorem} \ref{theor-G=l}).

\item For higher power gradings, $\ell_{n^{\prime}}(n^{\prime}-1)=\ell
_{n}(n-1)$ ({Theorem} \ref{theor-l=l}).
\end{itemize}

These ``quantization'' rules represent a distinctive feature of polyadic
grading, with no analog in binary graded algebra~theory.

3. \textit{{Support} properties}: We established that for strongly $\mathsf{G}%
$-graded polyadic algebras, the~support satisfies $\left\vert \text{supp}%
(\mathcal{A}^{[m,n]})\right\vert =|G|$, ensuring that all homogeneous
components   are~nontrivial.

4. \textit{{Graded} homomorphisms}: We developed the theory of polyadic graded
homomorphisms as pairs $(\Phi,\Psi)$, preserving both the algebraic structure and
grading, with~particular attention to cases where $\Psi=\operatorname{id}$.
The First Isomorphism Theorem for graded polyadic algebras is~proved.

\textbf{{Key} examples and applications}

1. \textit{{Ternary} superalgebras}: We constructed explicit examples including the following:

\begin{itemize}
\item Derived ternary superalgebras with binary $\mathbb{Z}_{2}$-grading;

\item Strictly nonderived ternary superalgebras graded by ternary groups
without identity ({Example} \ref{ex-tern}).
\end{itemize}

These examples demonstrate the possibility of grading by groups that lack
neutral elements---a situation impossible in classical grading~theory.

2. \textit{{Polynomial} algebras over }$n$\textit{-ary matrices}: We developed
the theory of polynomials over block-shift matrices, showing they can be
graded by polyadic integers $\mathbb{Z}^{[m^{\prime\prime},n^{\prime\prime}%
]}(a,b)$ with specific ``quantization'' conditions: $a=1$, $b=n-1$, and~
$k=\ell_{k}$ ({Section} \ref{sec-poly}).

3. \textit{{Higher} power gradings}: We introduced the concept of higher power
gradings where $n^{\prime}\neq n$, providing explicit solutions to the
``quantization'' condition $\ell_{n^{\prime}}(n^{\prime}-1)=\ell_{n}(n-1)$ and
constructing a $5$-ary superalgebra with ternary grading ({Section}
\ref{sec-high}).

\textbf{{Theoretical} implications}

1. \textit{Arity freedom with constraints}: While initial arities can be
chosen freely according to the arity freedom principle, meaningful grading
structures impose specific constraints through ``quantization'' rules. This
represents a natural selection mechanism for compatible arities in polyadic
graded~systems.

2. \textit{Role of identity elements}: The classical requirement that grading
groups contain an identity element $\mathsf{e}$ is relaxed in the polyadic
setting. We have shown that meaningful gradings exist even when $\mathsf{G}%
^{[n^{\prime}]}$ lacks such an element, as~demonstrated by strictly nonderived
ternary grading~groups.

3. \textit{Support as structural invariant}: The support of a graded polyadic
algebra becomes a more delicate invariant than in the binary case, intimately
connected with polyadic powers and word length ``quantization''.

4. \textit{Nonderived operations}: Our emphasis on strictly nonderived
operations (\textit{Remark} \ref{rem-nonder}) ensures that the constructed
examples are genuinely polyadic rather than disguised binary structures,
revealing the essential features of higher~arity.

\textbf{{Open} problems and future directions}

1. \textit{{Classification} problem}: Classify all possible multiary gradings
for given arities $(m,n)$ and finite group orders $|G|$, including the
enumeration of non-isomorphic \mbox{graded~structures.}

2. \textit{Homological aspects}: Develop homology and cohomology theories for
multiary graded algebras, extending the classical theory of group cohomology
for graded~rings.

3. \textit{Representation theory}: Investigate representations of multiary
graded algebras, particularly the decomposition of modules into homogeneous
components and the behavior of graded module~categories.

4. \textit{Physical applications}: Explore potential applications in
theoretical physics, where

\begin{itemize}
\item {Ternary and higher-arity structures appear in Nambu mechanics and
ternary ``quantization'';}

\item Higher power gradings might model symmetries in extended physical~systems;

\item Polyadic superalgebras could provide mathematical frameworks for
generalizations of supersymmetry.
\end{itemize}

5. \textit{Geometric realizations}: Develop geometric interpretations of
multiary graded algebras, potentially through noncommutative geometry or
graded manifold theory extended to polyadic~settings.

6. \textit{Connection with higher categories}: Investigate relationships
between multiary graded algebras and higher categorical structures, where
composition laws themselves have higher~arity.

7. \textit{Computational aspects}: {Develop algorithms for determining whether
a given polyadic algebra admits nontrivial multiary gradings, and~for
classifying such \mbox{gradings~computationally.}}

\textbf{{Concluding} remarks}

The theory of multiary graded polyadic algebras represents a natural and rich
extension of classical grading concepts to higher-arity algebraic structures.
The ``quantization'' rules we have discovered reveal an intricate interplay
between the arities of operations, the~order of grading groups, and~polyadic
powers---a layer of structure absent in   binary~algebra.

The examples we have constructed, particularly the strictly nonderived cases,
demonstrate that this generalization is not merely formal but yields genuinely
new mathematical objects with interesting properties. The~relaxation of
requirements such as the existence of identity elements in grading groups
opens new possibilities for algebraic~structures.

As polyadic algebra continues to develop as a field~\cite{duplij2022,dup2025h}%
, the~theory of graded structures will likely play an increasingly important
role, both in pure mathematics and in potential applications to physics and
other sciences. The~framework established here provides a solid foundation for
these future developments, offering both concrete examples and general
principles to guide further~exploration.

\pagestyle{emptyf}
\mbox{}


\begin{thebibliography}{}

\bibitem[\protect\citeauthoryear{Bahturin, Sehgal, and
  Zaicev}{\textcolor{blue}{\sc Bahturin et~al.}}{2001}]{bah/seh/zai}
{\textcolor{blue}{\sc Bahturin, Y., S.~Sehgal, and M.~Zaicev}} (2001).
\newblock Group gradings on associative algebras.
\newblock {\em J. Algebra\/}~{\bf 241}, 677--698.

\bibitem[\protect\citeauthoryear{Berezin}{\textcolor{blue}{\sc
  Berezin}}{1987}]{berezin}
{\textcolor{blue}{\sc Berezin, F.~A.}} (1987).
\newblock {\em Introduction to Superanalysis}.
\newblock Dordrecht: Reidel.

\bibitem[\protect\citeauthoryear{Bernstein, Leites, Molotkov, and
  Shander}{\textcolor{blue}{\sc Bernstein et~al.}}{2013}]{leites13}
{\textcolor{blue}{\sc Bernstein, J., D.~Leites, V.~Molotkov, and V.~Shander}}
  (2013).
\newblock {\em Seminars of Supersymmetries. Vol.1. Algebra and calculus}.
\newblock Moscow: MCCME.
\newblock In Russian, the English version is available for perusal.

\bibitem[\protect\citeauthoryear{Boboc, Dăscălescu, and {van
  Wyk}}{\textcolor{blue}{\sc Boboc et~al.}}{2024}]{bob/das/wyk}
{\textcolor{blue}{\sc Boboc, C., S.~Dăscălescu, and L.~{van Wyk}}} (2024).
\newblock Cyclic algebras, symbol algebras and gradings on matrices.
\newblock {\em Linear Algebra Appl.\/}~{\bf 688}, 157--178.

\bibitem[\protect\citeauthoryear{Bovdi}{\textcolor{blue}{\sc
  Bovdi}}{1974}]{bovdi}
{\textcolor{blue}{\sc Bovdi, A.~A.}} (1974).
\newblock {\em Group Rings}.
\newblock Uzhgorod: Uzgorod. Univ.

\bibitem[\protect\citeauthoryear{Cohn}{\textcolor{blue}{\sc Cohn}}{1965}]{cohn}
{\textcolor{blue}{\sc Cohn, P.~M.}} (1965).
\newblock {\em Universal Algebra}.
\newblock New York: Harper \& Row.

\bibitem[\protect\citeauthoryear{Deligne, Etingof, Freed, Jeffrey, Kazhdan,
  Morgan, Morrison, and Witten}{\textcolor{blue}{\sc Deligne
  et~al.}}{1999}]{del/eti/fre}
{\textcolor{blue}{\sc Deligne, P., P.~Etingof, D.~S. Freed, L.~C. Jeffrey,
  D.~Kazhdan, J.~W. Morgan, D.~R. Morrison, and E.~Witten}} (Eds.) (1999).
\newblock {\em Quantum Fields and Strings: A Cource for Mathematicians}, Vol.
  1, 2, Providence. American Mathematical Society.

\bibitem[\protect\citeauthoryear{D\"ornte}{\textcolor{blue}{\sc
  D\"ornte}}{1929}]{dor3}
{\textcolor{blue}{\sc D\"ornte, W.}} (1929).
\newblock Unterschungen \"uber einen verallgemeinerten {G}ruppenbegriff.
\newblock {\em Math. Z.\/}~{\bf 29}, 1--19.

\bibitem[\protect\citeauthoryear{Duplij}{\textcolor{blue}{\sc
  Duplij}}{2017}]{dup2017a}
{\textcolor{blue}{\sc Duplij, S.}} (2017).
\newblock Polyadic integer numbers and finite $(m,n)$-fields.
\newblock {\em p-Adic Numbers, Ultrametric Analysis and Appl.\/}~{\bf 9} (4),
  257--281.
\newblock {arXiv:math.RA/1707.00719}.

\bibitem[\protect\citeauthoryear{Duplij}{\textcolor{blue}{\sc
  Duplij}}{2019}]{dup2019}
{\textcolor{blue}{\sc Duplij, S.}} (2019).
\newblock Arity shape of polyadic algebraic structures.
\newblock {\em J. Math. Physics, Analysis, Geometry\/}~{\bf 15} (1), 3--56.

\bibitem[\protect\citeauthoryear{Duplij}{\textcolor{blue}{\sc
  Duplij}}{2022a}]{duplij2022}
{\textcolor{blue}{\sc Duplij, S.}} (2022a).
\newblock {\em Polyadic Algebraic Structures}.
\newblock London-Bristol: IOP Publishing.

\bibitem[\protect\citeauthoryear{Duplij}{\textcolor{blue}{\sc
  Duplij}}{2022b}]{dup2022a}
{\textcolor{blue}{\sc Duplij, S.}} (2022b).
\newblock Polyadization of algebraic structures.
\newblock {\em Symmetry\/}~{\bf 14} (9), 1782.

\bibitem[\protect\citeauthoryear{Duplij}{\textcolor{blue}{\sc
  Duplij}}{2025}]{dup2025h}
{\textcolor{blue}{\sc Duplij, S.}} (2025).
\newblock Higher power polyadic group rings, {\it preprint}  Univ. M\"unster,
  CIT,  M\"unster,  18~p., math.RA/2510.14029.

\bibitem[\protect\citeauthoryear{Green, Schwarz, and
  Witten}{\textcolor{blue}{\sc Green et~al.}}{1987}]{gre/sch/wit}
{\textcolor{blue}{\sc Green, M.~B., J.~H. Schwarz, and E.~Witten}} (1987).
\newblock {\em Superstring Theory}, Vol. 1,2.
\newblock Cambridge: Cambridge Univ. Press.

\bibitem[\protect\citeauthoryear{Hazewinkel and Gubareni}{\textcolor{blue}{\sc
  Hazewinkel and Gubareni}}{2016}]{haz/gub}
{\textcolor{blue}{\sc Hazewinkel, M. and N.~M. Gubareni}} (2016).
\newblock {\em Algebras, Rings and Modules: Non-commutative Algebras and
  Rings}.
\newblock Boca Raton: CRC Press.

\bibitem[\protect\citeauthoryear{Hazrat}{\textcolor{blue}{\sc
  Hazrat}}{2016}]{hazrat}
{\textcolor{blue}{\sc Hazrat, R.}} (2016).
\newblock {\em Graded Rings and Graded {G}rothendieck Groups}.
\newblock Cambridge: Cambridge University Press.

\bibitem[\protect\citeauthoryear{Kac}{\textcolor{blue}{\sc Kac}}{1977}]{kac3}
{\textcolor{blue}{\sc Kac, V.~G.}} (1977).
\newblock Lie superalgebras.
\newblock {\em Adv. Math.\/}~{\bf 26} (1), 8--96.

\bibitem[\protect\citeauthoryear{Kaku}{\textcolor{blue}{\sc
  Kaku}}{1998}]{kaku3}
{\textcolor{blue}{\sc Kaku, M.}} (1998).
\newblock {\em Introduction to Superstrings and $M$-Theory}.
\newblock Berlin: Springer-Verlag.

\bibitem[\protect\citeauthoryear{Kasner}{\textcolor{blue}{\sc
  Kasner}}{1904}]{kas}
{\textcolor{blue}{\sc Kasner, E.}} (1904).
\newblock An extension of the group concept.
\newblock {\em Bull. Amer. Math. Soc.\/}~{\bf 10}, 290--291.
\newblock Reported by L.~G.~Weld at 53rd Annual Meeting of AAAS, St. Louis,
  1904.

\bibitem[\protect\citeauthoryear{Kelarev}{\textcolor{blue}{\sc
  Kelarev}}{2002}]{kelarev}
{\textcolor{blue}{\sc Kelarev, A.}} (2002).
\newblock {\em Ring Constructions and Applications}.
\newblock Singapore: World Scientific.

\bibitem[\protect\citeauthoryear{Kumduang and
  Wattanatripop}{\textcolor{blue}{\sc Kumduang and
  Wattanatripop}}{2025}]{kum/wat}
{\textcolor{blue}{\sc Kumduang, T. and K.~Wattanatripop}} (2025).
\newblock Partial algebras of formulas under generalized superpositions.
\newblock {\em The Bulletin of Irkutsk State University. Series
  Mathematics\/}~{\bf 54}, 160--175.

\bibitem[\protect\citeauthoryear{Lam}{\textcolor{blue}{\sc
  Lam}}{1991}]{lam1991}
{\textcolor{blue}{\sc Lam, T.~Y.}} (1991).
\newblock {\em A First Course in Noncommutative Rings}.
\newblock New York: Springer.

\bibitem[\protect\citeauthoryear{Lang}{\textcolor{blue}{\sc Lang}}{2002}]{lang}
{\textcolor{blue}{\sc Lang, S.}} (2002).
\newblock {\em Algebra\/} (3rd ed.).
\newblock New York: Springer.

\bibitem[\protect\citeauthoryear{Leeson and Butson}{\textcolor{blue}{\sc Leeson
  and Butson}}{1980}]{lee/but}
{\textcolor{blue}{\sc Leeson, J.~J. and A.~T. Butson}} (1980).
\newblock On the general theory of $(m,n)$ rings.
\newblock {\em Algebra Univers.\/}~{\bf 11}, 42--76.

\bibitem[\protect\citeauthoryear{Matsumura}{\textcolor{blue}{\sc
  Matsumura}}{1959}]{mat59}
{\textcolor{blue}{\sc Matsumura, H.}} (1959).
\newblock On the graded rings of polynomials.
\newblock {\em J. Math. Soc. Japan\/}~{\bf 11} (1), 17--27.

\bibitem[\protect\citeauthoryear{N{\u{a}}st{\u{a}}sescu and
  Van~Oystaeyen}{\textcolor{blue}{\sc N{\u{a}}st{\u{a}}sescu and
  Van~Oystaeyen}}{1982}]{nas/oys}
{\textcolor{blue}{\sc N{\u{a}}st{\u{a}}sescu, C. and F.~Van~Oystaeyen}} (1982).
\newblock {\em Graded Ring Theory}.
\newblock Elsevier.

\bibitem[\protect\citeauthoryear{N{\u{a}}st{\u{a}}sescu and
  Van~Oystaeyen}{\textcolor{blue}{\sc N{\u{a}}st{\u{a}}sescu and
  Van~Oystaeyen}}{2004}]{nas/oys04}
{\textcolor{blue}{\sc N{\u{a}}st{\u{a}}sescu, C. and F.~Van~Oystaeyen}} (2004).
\newblock {\em Methods of Graded Rings}.
\newblock Berlin: Springer-Verlag.

\bibitem[\protect\citeauthoryear{Passman}{\textcolor{blue}{\sc
  Passman}}{1977}]{passman}
{\textcolor{blue}{\sc Passman, D.~S.}} (1977).
\newblock {\em The Algebraic Structure of Group Rings}.
\newblock New York-London-Sydney: John Wiley and Sons.

\bibitem[\protect\citeauthoryear{Post}{\textcolor{blue}{\sc Post}}{1940}]{pos}
{\textcolor{blue}{\sc Post, E.~L.}} (1940).
\newblock Polyadic groups.
\newblock {\em Trans. Amer. Math. Soc.\/}~{\bf 48}, 208--350.

\bibitem[\protect\citeauthoryear{Rotman}{\textcolor{blue}{\sc
  Rotman}}{2010}]{rotman}
{\textcolor{blue}{\sc Rotman, J.~J.}} (2010).
\newblock {\em Advanced Modern Algebra\/} (2nd ed.).
\newblock Providence: AMS.

\bibitem[\protect\citeauthoryear{Russell}{\textcolor{blue}{\sc
  Russell}}{1994}]{russell}
{\textcolor{blue}{\sc Russell, P.}} (1994).
\newblock Gradings of polynomial rings.
\newblock In: {\textcolor{blue}{\sc C.~L. Bajaj}} (Ed.), {\em Algebraic
  Geometry and its Applications}, New York: Springer, pp.\  365--373.

\bibitem[\protect\citeauthoryear{Sampson and George}{\textcolor{blue}{\sc
  Sampson and George}}{2026}]{sam/geo}
{\textcolor{blue}{\sc Sampson, M. and R.~George}} (2026).
\newblock Semigroups in distributed computations: n-ary operations and
  irreducibility.
\newblock {\em Eur. J. Pure Appl. Math.\/}~{\bf 19} (1), 7215.

\bibitem[\protect\citeauthoryear{Sokhatsky}{\textcolor{blue}{\sc
  Sokhatsky}}{1997}]{sok1}
{\textcolor{blue}{\sc Sokhatsky, F.~M.}} (1997).
\newblock On the associativity of multiplace operations.
\newblock {\em Quasigroups Relat. Syst.\/}~{\bf 4}, 51--66.

\bibitem[\protect\citeauthoryear{Terning}{\textcolor{blue}{\sc
  Terning}}{2005}]{terning}
{\textcolor{blue}{\sc Terning, J.}} (2005).
\newblock {\em Modern Supersymmetry}.
\newblock Oxford: Oxford University Press.

\bibitem[\protect\citeauthoryear{Wess and Bagger}{\textcolor{blue}{\sc Wess and
  Bagger}}{1983}]{wes/bag}
{\textcolor{blue}{\sc Wess, J. and J.~Bagger}} (1983).
\newblock {\em Supersymmetry and Supergravity}.
\newblock Princeton: Princeton Univ. Press.

\bibitem[\protect\citeauthoryear{Zalesskii and Mikhalev}{\textcolor{blue}{\sc
  Zalesskii and Mikhalev}}{1975}]{zal/mik}
{\textcolor{blue}{\sc Zalesskii, A.~E. and A.~V. Mikhalev}} (1975).
\newblock Group rings.
\newblock {\em J. Soviet Math.\/}~{\bf 4} (1), 1--78.

\end{thebibliography}
\end{document}